\documentclass[12pt,leqno]{amsart}

\usepackage{amssymb}
\usepackage{amsthm}
\usepackage{amsmath}
\usepackage{amscd}
\usepackage[mathscr]{eucal}

\numberwithin{equation}{section}

\theoremstyle{plain}
\newtheorem{theorem}{Theorem}[section]
\newtheorem{corollary}[theorem]{Corollary}
\newtheorem{lemma}[theorem]{Lemma}
\newtheorem{proposition}[theorem]{Proposition}

\theoremstyle{definition}
\newtheorem{definition}[theorem]{Definition}
\newtheorem{remark}[theorem]{Remark}

\theoremstyle{remark}

\newcommand{\R}{\mathbb{R}}
\newcommand{\Q}{\mathbb{Q}}
\newcommand{\Z}{\mathbb{Z}}

\newcommand{\C}{\mathbb{C}}
\newcommand{\h}{\mathbb{H}}

\newcommand{\G}{\Gamma}
\newcommand{\g}{\gamma}

\newcommand{\La}{\Lambda}

\newcommand{\back}{\backslash}

\newcommand{\V}{V(\Q)}

\begin{document}

\title[Cycles in Hyperbolic Manifolds and Siegel Modular Forms]
       {Cycles in Hyperbolic Manifolds of Non-compact
       Type and Fourier Coefficients of Siegel Modular Forms}

\date{}

\author[Jens Funke]{Jens Funke*}
\address{Department of Mathematics \\
    Rawles Hall \\
    Indiana University \\
    Bloomington, IN 47405 \\
    USA}
\email{jefunke@indiana.edu}
\thanks{*Fellow of the European PostDoc Institute (EPDI) 99-01}

\author[John Millson]{John Millson**}
\address{ Department of Mathematics \\
   University of Maryland \\
   College Park, MD 20742 \\
   USA}
\email{jjm@math.umd.edu}
\thanks{**Partially supported by NSF-grant DMS 98-03520}

\maketitle

\section{Introduction}

Throughout the 1980's, Kudla and the second named author studied
integral transforms $\Lambda$ from closed differential forms on
arithmetic quotients of the symmetric spaces of orthogonal and
unitary groups to spaces of classical Siegel and Hermitian modular
forms (\cite{KMI,KMII,KMCan,KM90}). These transforms came from the
theory of dual reductive pairs and the theta correspondence.

In \cite{KM90} they computed the Fourier expansion of
$\Lambda(\eta)$ in terms of periods of $\eta$ over certain totally
geodesic cycles under the assumption that $\eta$ was rapidly
decreasing. This also gave rise to the realization of intersection
numbers of these `special' cycles with cycles with compact support
as Fourier coefficients of modular forms.

It is clear from \cite{HZ},\cite{Cogdell} and \cite{Funke} that
the situation is far more complicated when the hypothesis of rapid
decay is dropped. The purpose of this paper is to initiate a
systematic study of this transform for non rapidly decreasing
differential forms $\eta$ by considering  the case for the finite
volume quotients of hyperbolic space coming from unit groups of
isotropic quadratic forms over $\Q$. We expect that many of the
techniques and features of this case will carry over to the more
general situation.

\vspace{.5cm}

We now give a more precise description of this paper. Let $V(\Q)$
be a rational vector space of dimension $m=p+1$ with a symmetric
bilinear form $(\;,\;)$ of signature $(p,1)$ and put $G(\Q) =
SO(V(\Q))$. We let $L$ be an integral lattice in $V(\Q)$ and
$\G(\Q)$ be a torsion-free subgroup of the stabilizer of $L$ in
$G(\Q)$. We denote by $B$ the associated symmetric space to
$G(\R)$, and we assume that the hyperbolic manifold $M=\G \back B$
is non-compact.

Kudla and the second named author (\cite{KMI,KMII}) constructed a
certain theta function $\theta(\tau,Z)$ for $\tau \in \h_n$, the
Siegel upper half space, and $Z \in B$, which is a
\emph{non-holomorphic} Siegel modular form of weight
$\tfrac{m}{2}$
with values in the closed differential $n$-forms of $M$. For
$\eta$ a rapidly decreasing closed differential $(p-n)$-form in
$M$, they then defined the transform
\begin{equation}\label{transform}
\Lambda(\eta)(\tau) = \int_{M} \eta \wedge \theta(\tau,Z).
\end{equation}
They showed that $\Lambda(\eta)(\tau)$ is a holomorphic cusp form,
see \cite{KM90}. Moreover, the Fourier coefficients are given as
periods of $\eta$ over certain geometrically defined composite, in
general non-compact, `special' cycles $C_{\beta}$  in $M$ attached
to positive definite $\beta \in Sym_n(\Q)$, i.e.,
\begin{equation}
 \Lambda(\eta)(\tau) = \sum_{\beta > 0} \left(  \int_{C_{\beta}} \eta
 \right) e^{2\pi i tr(\beta \tau)}.
\end{equation}
The lift factors through the cohomology $H_c^{p-n}(M,\C)$ with
compact support, and the period $\int_{C_{\beta}} \eta$ is the
evaluation of the pairing of $[\eta] \in H_c^{p-n}(M,\C)$ with the
relative cycle $C_{\beta} \in H_{p-n}(M,\partial{M},\Z)$. The key
point is here that the Fourier coefficients $\theta_{\beta}$ of
$\theta(\tau)$ are the Poincar\'e-dual forms of the cycles
$C_{\beta}$.

In the case of $p=2$ and $n=1$, this lift is closely related to
the work of Shintani \cite{Shintani} on the inverse of the Shimura
lift.

\vspace{.5cm}

The Borel-Serre compactification makes $M$ a compact manifold with
boundary $\overline{M}$. Here each boundary component is a
$(p-1)$-torus at the various cusps of $M$. We develop a machinery
to determine the growth of $\Theta(\tau,Z)$ and show
\begin{theorem}\label{THA}\hfill

$\theta(\tau,Z)$ extends to a smooth differential form on
$\overline{M}$. Moreover, the coefficients of the restriction of
$\theta(\tau,Z)$ to each boundary component are given by a linear
combination of holomorphic Siegel cusp forms of weight
$\tfrac{m}{2}$ coming from the orthogonal group $O(p-1)$.

\end{theorem}

We can therefore extend the theta integral (\ref{transform}) to
$(p-n)$-forms $\eta$ on $\overline{M}$. For the special case $n=p$
and $\eta=1$, the theta integral was already studied by Kudla
(\cite{KudlaHS,KShin}).
\begin{theorem} \label{THB}\hfill

Let $\eta$ be a closed differential $(p-n)$-form on
$\overline{M}$. Then  $\Lambda(\eta)(\tau)$ is a
\underline{\emph{holomorphic}} Siegel modular form of weight
$\tfrac{m}2$ for a suitable congruence subgroup of $Sp(n,\Z)$.
\end{theorem}
The key point is here that there exists another, \emph{rapidly
decreasing}  theta function $\Xi(\tau,Z)$ such that
\begin{equation}
\bar{\partial} \, \theta(\tau,Z) = d \, \Xi(\tau,Z).
\end{equation}
Here $\bar{\partial}$ operates on the $\tau$-variable and $d$ on
the $Z$-variable. This, together with Stokes' theorem, implies
that $\Lambda(\eta)(\tau)$ satisfies the Cauchy-Riemann equations.

The form $\Xi$ exists in general but it is not necessarily rapidly
decreasing. Thus the problem of when $\Lambda(\eta)$ is
holomorphic is rather delicate. In fact, in \cite{Funke} it was
shown that in the case of signature $(p,2)$ analogous theta
integrals are in general non-holomorphic modular forms.

We call the space of holomorphic Siegel cusp forms of weight
$\tfrac{m}2$ and degree $n$ coming from theta series attached to
$O(p-1)$ the space of unstable cusp forms and denote it by
$\Theta^{(n)}(p-1)$. (For $p=2$ and $n=1$, these cusp forms
correspond to Eisenstein series of weight $2$ under the Shimura
correspondence).

By Theorem \ref{THA} the image of exact forms lies in the space of
unstable cusp forms. Denoting the space of holomorphic Siegel
modular forms of weight $\tfrac{m}2$ and degree $n$ by
$M^{(n)}_{m/2}$, we therefore obtain
\begin{theorem}\label{THC}\hfill

The transform $\Lambda$ factors through the cohomology
$H^{p-n}(\overline{M},\C) \simeq H^{p-n}(M,\C)$ modulo unstable
Siegel cusp forms, i.e., $\Lambda$ defines a map
\begin{equation*}
\Lambda: H^{p-n}(\overline{M},\C) \longrightarrow M^{(n)}_{m/2} /
\Theta^{(n)}(p-1).
\end{equation*}

\end{theorem}

By Theorem \ref{THB} we see by the Koecher principle that the
Fourier expansion of $\Lambda(\eta)(\tau)$ is given by
\begin{equation}
 \Lambda(\eta)(\tau) = \sum_{\beta \geq 0}  a_{\beta}(\eta)
e^{2\pi i tr(\beta\tau)}
\end{equation}
with
\begin{equation}\label{Thomintegral}
a_{\beta}(\eta) = \int_M \eta \wedge \theta_{\beta}(\tau).
\end{equation}
(For $n=1$, the vanishing of the negative coefficients follows
from a direct calculation which we omit).

For the singular coefficients, the $\theta_{\beta}(\tau)$ turn out
to be rapidly decreasing, and we have
\begin{theorem}\label{THD}
\begin{equation*}
a_{\beta}(\eta) =
\begin{cases}
0 \qquad \qquad \qquad  \; \; \text{if} \qquad rk(\beta) < n-1 \\
(-1)^n \int_{C^s_{\beta}} \eta \qquad \; \text{if} \qquad
rk(\beta) = n-1.
\end{cases}
\end{equation*}
\end{theorem}

In particular, we see that $\Lambda(\eta)(\tau)$ is in general no
longer a cusp form. Here, for $\beta$ positive semi-definite of
rank  $n-1$, the `singular' cycles $C^s_{\beta}$ are linear
combinations of embedded $(p-n)$-subtori at each component of the
Borel-Serre boundary of $M$. The coefficients are values of
Dirichlet series attached to the boundary components. Note that
the $C^s_{\beta}$ can be considered as absolute cycles in $M$ and
therefore the period of $\eta$ over  $C^s_{\beta}$ is
cohomological.

The calculation of the singular Fourier coefficients uses
extensively ideas from \cite{KShin}, where the case of $n=p$ was
considered. However, through a careful growth analysis of the
theta series involved we are able to greatly simplify the concept
of the calculations, avoiding the usage of a wave packet attached
to Eisenstein series. This observation should also be very helpful
for extending the much more general results of \cite{KM90}.

\vspace{.5cm}

The situation for the positive definite coefficients is
considerably more complicated as now $\theta_{\beta}$ is nonzero
at the boundary and therefore homotopy- and Stokes-type arguments
for the computation of (\ref{Thomintegral}) are no longer
available. In particular, the calculation for $\eta$ rapidly
decreasing (see \cite{KMCan}) does not extend to arbitrary $\eta$.
This corresponds to the fact that the period $ \int_{C_{\beta}}
\eta$ (where $C_{\beta}$ is the (in general relative) cycle
mentioned above) no longer has a (co)homological interpretation.

In fact, if $\eta$ is an exact form which extends to the boundary,
the equation
\begin{equation}\label{ThomL}
a_{\beta}(\eta) \stackrel{?}{=}  \int_{C_{\beta}} \eta
\end{equation}
is in general no longer valid! We define the `defect' 
$\delta_{\beta}(\eta) = a_{\beta}(\eta) - \int_{C_{\beta}} \eta$ and show that $\delta_{\beta}$  descends to a function on ${Z}^{p-n}(\overline{M})/ 
{Z}^{p-n}(\overline{M}, \partial \overline{M})$, where 
 ${Z}^{\ast}(\overline{M})$ is the space of closed differential forms on $\overline{M}$ and ${Z}^{\ast}(\overline{M}, \partial \overline{M})$ the subspace of forms which vanish at the boundary. Moreover, we show that the defect can be non-zero on the subspace of exact $(p-n)$-forms supported near 
 $\partial \overline{M}$.


\vspace{.5cm}

For the case of a Riemann surface, i.e., for the
case of $SO(2,1)$ and $n=1$, we have a complete picture:

\begin{theorem}\label{THE}
Let $p=2$ and $n=1$. Then each class in $H^1(\overline{M},\C)$ has
a representative $\eta$ such that
\begin{equation}\label{ThomEis}
\Lambda(\eta)(\tau)= \left(\int_{C^s_{0}} \eta \right)\; + \;
\sum_{\beta > 0} \left( \int_{C_{\beta}} \eta \right) e^{2\pi i
\beta \tau}.
\end{equation}
Hence (\ref{ThomEis}) holds in $M^{(1)}_{3/2} / \Theta^{(1)}(1)$
for all closed $1$-forms $\eta$ in $\overline{M}$.
\end{theorem}

The point is here that via the theory of Eisenstein cohomology
$H^1(\overline{M},\C)$ splits into its cuspidal (or $L_2$)
cohomology and a part defined by Eisenstein series coming from
cohomology classes at the boundary. We are able to directly
compute (\ref{Thomintegral}) for forms defined by cusp forms and
Eisenstein series, thus verifying (\ref{ThomL}).

Furthermore, we can consider the 'truncated' part $\theta^c(\tau)$
of the form $\theta(\tau)$, which is obtained by subtracting the
Eisenstein form of the restriction of $\theta(\tau)$ to the
boundary from $\theta(\tau)$ itself. $\theta^c(\tau)$ is again a
modular form of weight $3/2$ with values now in the rapidly
decreasing differential $1$-forms of the Riemann surface $M$.

For $\beta >0$, we  define $C^c_{\beta}$ to be the homology class
dual to the $\beta$-th Fourier coefficient of $\theta^c(\tau)$.
This definition and the following result is completely analogous to
the one by Hirzebruch-Zagier for Hilbert modular surfaces
(\cite{HZ}):

\begin{theorem}\label{THF}
Let $p=2$ and $n=1$. The map
\[
\eta \mapsto \int_M \eta \wedge \theta^c(\tau)
\]
factors through $H^1(\overline{M},\C)$, and if
 $C$ is the homology class dual to
$[\eta]$, we have that
\begin{equation*}
\int_M \eta \wedge \theta^c(\tau) = -[C^s_{0}.C] + \sum_{\beta >0}
[C^c_{\beta}.C] e^{2\pi i \beta \tau}
\end{equation*}
is a holomorphic modular form of weight $3/2$. Here $[\, . \,]$
denotes the cohomological intersection product.
\end{theorem}

It seems natural to expect that this generalizes to $SO(p,1)$ (at
least when the Eisenstein classes involved are not residual), and
we hope to come back to this issue in the near future.

\vspace{.5cm}

We can also define in the general case
\begin{equation}
\Lambda(C)(\tau) = \int_C \theta(\tau)
\end{equation}
for $C$ being a special cycle of complementary dimension $n$. For
this lift, we have complete control over the Fourier
coefficients:
\begin{theorem}\label{THG}
$\Lambda(C)(\tau)$ is a holomorphic Siegel modular form of weight
$\tfrac{m}2$ and degree $n$ and
\begin{equation*}
\Lambda(C)(\tau) = \sum_{\beta >0} [C.C_{\beta}]_{tr} e^{2\pi i
tr(\beta\tau)} \; + \; (-1)^n \sum_{\substack{\beta \geq 0 \\
rk(\beta) = n-1}} [C.C^s_{\beta}] e^{2\pi i tr(\beta\tau)}.
\end{equation*}
 Here  $[C.C_{\beta}]_{tr}$ denotes the transversal intersection
 number of $C$ and $C_{\beta}$ in $M$, i.e., the sum of the
 transversal intersections counted with multiplicities.
\end{theorem}

\vspace{.5cm}

We would like to thank Steve Kudla for many crucial discussions and
his encouragement. The first named author would like to thank the
Max-Planck-Institut f\"ur Mathematik in Bonn and the Department of
Algebra and Geometry at the University of Barcelona  for their
hospitality where major work for this paper was done.

\section{Preliminaries}

Let $V(\Q)$ be a rational vector space of dimension $m= p+1$ and
let $(\;,\;)$ be a  non-degenerate symmetric bilinear form on
$V(\Q)$ with signature $(p,1)$. Let $L \subset \V$ be an integral
$\Z$-lattice of full rank, i.e., $L \subset L^{\#}$, the dual
lattice. We let $G(\Q) = SO(\V)$ viewed as an algebraic group over
$\Q$. We denote by $\G(L)$ the stabilizer of the lattice $L$ and
fix a neat  subgroup $\G$ of finite index in $\G(L) \cap G_0(\R)$,
which acts trivially on $L^{\#}/L$. Here  $G_0(\R)$ is the
connected component of the identity of $G(\R)$.

Let $B$ be the real hyperbolic space of dimension $p$ and realize
$B$ as one component of the two-sheeted hyperboloid of vectors of
length $-1$:
\begin{equation}
B = \{Z \in V(\R):  (Z,Z) = -1\}^0.
\end{equation}

Fix a base point $Z_0 \in B$ and let $K$ be the stabilizer of
$Z_0$ in $G_0(\R)$. Then $K \simeq SO(p)$ is a maximal compact
subgroup of $G_0(\R)$, and we have
\begin{equation}
B \simeq G_0(\R)/K.
\end{equation}
Note that we can identify $B$ as the set of negative lines in
$V(\R)$ and therefore also as the space of minimal majorants of
$(\;,\;)$ by defining, for $Z \in B$, the majorant
\begin{equation}
(\;,\;)_Z =
\begin{cases}
(\;,\;) \; &\text{on $Z^{\perp}$;} \\
-(\;,\;) &\text{ on $\R \, Z$.}
\end{cases}
\end{equation}

For the tangent space $T_Z(B)$ we have the standard canonical
identification
\begin{equation}
T_Z(B) \simeq Z^{\perp}.
\end{equation}
We fix an orientation on $V$, and this induces an orientation of
$B$ by requiring that, for every properly oriented basis
$\{w_1,...,w_p\}$ for $T_Z(B) \simeq Z^{\perp}$, the basis
$\{w_1,...,w_p,Z\}$ is properly oriented for $V$. Note that the
action of $G_0(\R)$ on $B$ preserves this orientation.

\vspace{.5cm}

We assume that the hyperbolic manifold $M= \G \back B$ is
non-compact. It is well known \cite{Borel} that this is the case if
and only if $\V$ has an isotropic vector. Then $\G$ acts with
finitely many orbits on the set of isotropic lines in $\V$, the
cusps of $M$.
 We choose cusp representatives $\ell_0, \ell_1,...,\ell_r$ and primitive vectors $u_j \in L$ such that
\begin{equation}\label{forward}
\ell_j = \Q \, u_j \qquad \text{and} \qquad (u_j,Z) < 0
\end{equation}
for all $Z \in B$. We will express this second condition by saying
$u_i$ is forward pointing. We note that every null line has a
canonical orientation given by the class of a forward pointing
vector. We also choose $g_j \in G_0(\Q) = G(\Q) \cap G_0(\R)$ such
that
\begin{equation}
g_j \, u_0 = u_j
\end{equation}
and with $g_0 = 1$. Pick another isotropic vector $u'_0 \in \V$
such that $(u_0,u'_0) = -1/2$. This gives an isomorphism
\begin{equation}
\ell_0^{\perp}/\ell_0 \simeq W(\Q) := [u_0,u'_0]^{\perp}
\end{equation}
Note that $W$ is positive definite of dimension $p-1$. We choose a
basis $\{w_1,...,w_{p-1}\}$ of $W$ such that
$u_0,w_1,...,w_{p-1},u'_0$ is a positively oriented basis for $\V$
and call such a basis a Witt basis for $V(\Q)$. Note that this
also gives rise to an orientation of $\ell_0^{\perp} / \ell_0$.
With respect to this basis, $(\;,\;)$ is of the form
\begin{equation}
(\;,\;) \sim
\begin{pmatrix}
 &  & -1/2 \\
 & S&    \\
-1/2& &    \\
\end{pmatrix},
\end{equation}
where $S$ is the matrix of the bilinear form restricted to $W$.

We can assume that the base point $Z_0$ is rational and contained
in the hyperbolic plane $[u_0,u'_0]$. Since we assumed $(Z,Z) =
-1$ and $(Z,u_0) < 0$, we see that $Z_0 = u_0 + u'_0$, i.e., in
coordinates:
\begin{equation}
Z_0 =
\begin{pmatrix}
1\\0\\1
\end{pmatrix}.
\end{equation}
Note that majorant $(\;,\;)_{Z_0} =: (\;,\;)_0$ associated to the
base point $Z_0$ is given by
\begin{equation}
(\;,\;)_0 \sim
\begin{pmatrix}
1/2& & \\
 &S& \\
 & &1/2
\end{pmatrix}.
\end{equation}

We pick another basis for $V(\R)$ as follows. We let
\begin{equation}
e_1 = u_0 - u'_0 \qquad \qquad \text{and} \qquad \qquad e_{p+1} =
u_0 + u'_0 = Z_0.
\end{equation}
We have $(e_1,e_1) = 1$ and $e_1 \perp Z_0$ and extend $e_1$ to an
orthonormal basis $\{ e_1, \cdots, e_p\}$ for $Z_0^{\perp}$. With
respect to this basis $\{e_1, \cdots, e_{p+1} \}$ the bilinear
form has the matrix
\begin{equation}
(\;,\;) \sim
\begin{pmatrix}
1 &      &  &  \\
  &\ddots&  &  \\
  &      &1 &  \\
  &      &  &-1
\end{pmatrix}.
\end{equation}

Let $\mathfrak{g}$ be the Lie algebra of $G_0(\R)$ and
$\mathfrak{k}$ be that of $K$. We then have the Cartan
decomposition
\begin{equation}
\mathfrak{g} = \mathfrak{k} + \mathfrak{p},
\end{equation}
where $\mathfrak{p}$ is the orthogonal complement of
$\mathfrak{k}$ with respect to the Killing form. We identify
$\mathfrak{p}$ with $Z_0^{\perp}$ in an $SO(p)$-equivariant way
via
\begin{equation}
\begin{matrix}
Z_0^{\perp} &\xrightarrow{\sim} &\mathfrak{p} \\
v &\longmapsto & v \wedge Z_0,
\end{matrix}
\end{equation}

where $w \wedge w' \in \bigwedge^2 V$ is identified with an
element of $\mathfrak{g}$ given by
\begin{equation}
(w \wedge w')(v) = (w,v)w' - (w',v)w.
\end{equation}

 We identify the basis $\{e_1, \cdots, e_p\}$ for $Z_0^{\perp}$ with a basis of $\mathfrak{p}$. With respect to this basis we have
\begin{equation}
\mathfrak{p} \simeq \left\{ \begin{pmatrix} 0&v \\ ^tv&0
\end{pmatrix}: v \in Z_0^{\perp} \right\}.
\end{equation}
We let $\{\omega_1,\cdots,\omega_p\}$ be the dual basis of
$\mathfrak{p}^{\ast}$ corresponding to this basis.

We will denote coordinates with respect to the Witt basis
$\{u_0,w_1,...,w_{p-1},u'_0\}$ with $y_{ij}$ and coordinates with
respect to the basis $\{e_i\}$ with $x_{ij}$.

\vspace{.5cm}

Let $P$ be the $\Q\;$-parabolic subgroup of $G$ defined by
\begin{equation}
P(\Q) = \{ g \in G(\Q) : g\ell_0 = \ell_0 \}.
\end{equation}
Then for the unipotent radical $N(\Q)$ of $P(\Q)$, we have
\begin{equation}
N(\Q) \simeq W(\Q),
\end{equation}
and the isomorphism is explicitly given by
\begin{equation}
N(\Q) \simeq
\begin{matrix} \left\{ n(w)  =
\begin{pmatrix}
1&2(\cdot,w)& (w,w) \\
 &   1_W   & w  \\
 &         & 1  \\
\end{pmatrix}
: \; w \in W(\Q) \right\}
\end{matrix}.
\end{equation}

The maximal $\Q\;$-split torus $A(\Q)$ is given by
\begin{equation}
A(\Q) \simeq \left\{ a(t) =
\begin{pmatrix}
t& &  \\
 &   1_W   &  \\
 &         & t^{-1}  \\
\end{pmatrix}
: \; t \in \Q \right\}.
\end{equation}

We define
\begin{equation}
M = P(\R) \cap K \simeq SO(W(\R))
\end{equation}
and have the standard decompositions
\begin{equation}
G_0(\R) = N(\R) \, A_0(\R) \, K
\end{equation}
and
\begin{equation}
P_0(\R) = N(\R) \, A_0(\R) \, M,
\end{equation}
where $P_0(\R) = P(\R) \cap G_0(\R)$ and $A_0(\R) = A(\R) \cap
G_0(\R) \simeq \R_+$.

For $t \in \R_+$, let
\begin{equation}
A_t = \{ a(t') \in A_0(\R) : t' > t \},
\end{equation}
and for an open relatively compact subset $\omega \subset N(\R)$,
define the Siegel set
\begin{equation}
\mathfrak{S}_t = \omega A_t K \subset G_0(\R).
\end{equation}
Then by \cite{Borel} there exists a Siegel set $\mathfrak{S}
\subset G_0({\R})$ such that
\begin{equation}
G_0(\R) = \bigcup_{j} \G g_j \mathfrak{S}
\end{equation}
and
\begin{equation}
B= \bigcup_{j} \G g_j \mathfrak{S}',
\end{equation}
where $ \mathfrak{S}' = \mathfrak{S} \, \cdot \, Z_0.$

Let $N_j$, $0 \leq j \leq r$, be the point-wise stabilizer of the
cusps $\ell_j = \Q \, u_j$ in $N$ and $\G_j = N_j \cap \G$. We
have
\begin{equation}
N_j = g_j N g_j^{-1}.
\end{equation}
There exist lattices $\Lambda_j \subset W(\Q)$ such that
\begin{equation}\label{lattices}
\G_j = \{ g_j n(\lambda) g_j^{-1} : \lambda \in \Lambda_j \}.
\end{equation}

Recall that by adding for each cusp $\ell_j$ the torus $\G_j \back
N_j$ to the manifold $ M= \G \back B$ we obtain (with the
appropriate topology) the compact manifold with boundary $\bar{
M}$. This is the Borel-Serre compactification, see \cite{BS}. We
have
\begin{align}
    \overline{M} &= M \coprod_{j=0}^r \G_j \back N_j.
\end{align}

\vspace{.5cm}

We introduce upper-half space coordinates on $B$ associated to an
isotropic line, which we take to be $\ell_0$. We consider the map
\begin{equation}
\sigma: A_0(\R) \times N(\R) \longrightarrow B
\end{equation}
given by
\begin{equation}
\sigma(a,n) = n \, a \, Z_0
\end{equation}
Via the parametrization of $ A_0(\R) \times N(\R)$ by $\R_+ \times
\R^{p-1}$ we obtain coordinates on $B$ by
\begin{equation}
(t,b) \longmapsto Z(t,b) := n(b) \, a({t}) \, Z_0.
\end{equation}
We have
\begin{equation}
Z(t,b) =
\begin{pmatrix}
{t} +  t^{-1} (b,b) \\
t^{-1} b \\
t^{-1}
\end{pmatrix},
\end{equation}
where we identified $\R^{p-1}$ with $W(\R) \simeq N(\R)$. We
observe that in $\mathbb{P}(V)$ we have
\begin{equation}
     \lim_{t \to \infty} Z(t,b) = \ell_0,
\end{equation}
whereas in the Borel-Serre enlargement of $B$ we have $\lim_{t \to
\infty } Z(t,b) =b \in \ell_0^{\perp} / \ell_0$.

We extend $\sigma$ to $ N \times A \times K \longrightarrow G$
 by $\sigma (n,a,k) = nak$, and this induces an isomorphism between the left-invariant forms on $NA$ and the horizontal left-invariant forms on $G$ which we identify with $\mathfrak{p}^{\ast}$. It is easily seen that a basis for the left-invariant forms on $NA$ is given (in terms of the left-invariant forms  $\frac{dt}{t}$ on $A$ and $db_i$, $ 1 \leq i \leq p-1$ on $N$) by $\{\nu_1,\nu_2,\cdots,\nu_p\}$, where
\begin{equation}
  \nu_1 = \frac{dt}{t} \qquad \qquad \text{and} \qquad \qquad \nu_i = \frac{db_  {i-1}}{t} \qquad \text{for} \qquad 2 \leq i \leq p.
\end{equation}
We have
\begin{lemma}\label{invariantforms}
\[
\sigma^{\ast} \; \omega_i = \nu_i \qquad \qquad \text{for} \quad 1
\leq i \leq p.
\]
\end{lemma}

\begin{proof}
We only have to prove this at the identity. Then the basis
$\{\nu_1,\cdots,\nu_p\}$ for $\mathfrak{a}^{\ast} \, + \,
\mathfrak{n}^{\ast}$ is dual to the basis $\{ 2 u_0 \wedge u_0', 2
e_j \wedge u_0, 2\leq j \leq p \}$. The image of this basis under
$d\sigma |_e$ when projected onto $\mathfrak{p}$ (the horizontal
Maurer-Cartan forms annihilate $\mathfrak{k}$) is the basis $\{
e_1 \wedge e_{p+1},  e_2 \wedge e_{p+1}, \cdots,  e_p \wedge
e_{p+1} \}$. But this basis is dual to $\{\omega_i\}$ per
definitionem.
\end{proof}

\vspace{.5cm}

We will need a refinement of these coordinates associated to
positive semi-definite subspace $U$ of $V(\Q)$ of dimension $n$
such that the radical
\begin{equation}
R(U) = \{ u \in U : (u,U) =0 \}
\end{equation}
is non-zero. In this case we see by signature considerations that
there exists a rational isotropic line $\ell$ such that
\begin{equation}
R(U) = \ell.
\end{equation}
We may choose the above Witt decomposition such that
\begin{equation}
U = \ell_0 \; + U \cap W,
\end{equation}
i.e., $\ell = \ell_0$. We write $U' = U \cap W$ and let $U''$ be
the orthogonal complement of $U'$ in $W$, hence
\begin{equation}
W = U' \oplus U'',
\end{equation}
with the summands orthogonal for both $(\;,\;)$ and $(\;,\;)_0$.

We define subgroups $N'$ and $N''$ of $N$ with
\begin{equation}
N' \simeq U' \qquad \text{and} \qquad N'' \simeq U''
\end{equation}
under the isomorphism from $W$ to $N$. We also define
\begin{align*}
N_U &= \{ n \in N : n|_U = id \} = \{ n \in N : n |_{U'} = id \}.
\end{align*}
We observe that
\begin{equation}
N_U = N''.
\end{equation}
Indeed, for $w,w' \in W$, we have $n(w)w' = w' + (w,w') u_0$,
whence $N_U= (U')^{\perp} = U''$.

We can write
\begin{equation}
n(b) = n(b') n(b'')
\end{equation}
with $n(b') \in N'$ and $n(b'') \in N''$; so $b' \in \R^{n-1}
\simeq U'$ and $b'' \in \R^{p-n} \simeq U''.$ We obtain a product
decomposition
\begin{equation}
\sigma : \R_+ \times \R^{n-1} \times \R^{p-n} \longrightarrow B
\end{equation}
with
\begin{equation}
\sigma(t,b',b'') = Z(t,b',b'') := n(b')n(b'')a(t) Z_0.
\end{equation}

\section{Special Cycles}\label{cycles}

We define special cycles in $B$ as follows: Let $U$ be a positive
definite subspace of $V(\R)$ of dimension $n \leq p$, and define
\begin{equation}
B_U = \{ Z \in B : Z \perp U\}.
\end{equation}
Note that $B_U$ is a totally geodesic  submanifold, isomorphic to
the hyperbolic space of dimension $p-n$. If $U = span_{\R}X$ for
an $n$-frame $X=(x_1,\cdots,x_n)$ in $V(\R)$, we also write $B_X$
for $B_U$. An orientation on $U$ (say, coming from $X$) induces
one on $B_U$ as follows. We have a canonical isomorphism
\begin{equation}
T_Z(B_U) \simeq Z^{\perp} \cap U^{\perp}.
\end{equation}
Then $T_Z(B_U)$ receives an orientation by the rule that the
orientation of $T_Z(B_U)$ followed by the orientation of $U = U
\cap Z^{\perp}$ is the orientation of $T_Z(B) \simeq Z^{\perp}$.

Let $G_U$ be the point-wise stabilizer of $U$ in $G$ and put $\G_U
= \G \cap G_U$. We then define $C_U = \G_U \backslash B_U$; the
image of $B_U$ in $M$.

 For $\beta \in Sym_n(\Q)$, we consider the corresponding hyperboloid
\begin{equation}
\Omega_{\beta} = \{ X \in V(\Q) : \frac12 (X,X) = \beta \},
\end{equation}
with $(X,X)_{ij} = (x_i,x_j)$.

We fix a congruence condition $h \in \left(L^{\#}\right)^n$ once
and for all.

If $\beta$ is positive definite, then $\G$ acts on $\Omega_{\beta}
\cap (L^n +h) $ with finitely many orbits, and we define the
composite cycle
\begin{equation}
C_{\beta} = \sum_{\G \backslash \Omega_{\beta} \cap (L^n +h)} C_X.
\end{equation}

\vspace{.5cm}

We now construct special cycles on the Borel-Serre boundary of
$\overline{M}$. Let $U$ be a positive semidefinite subspace of
$V(\Q)$ of dimension $n$ with nonzero radical $R(U) = \ell$. We
denote the unipotent radical of the parabolic associated to $\ell$
by $N_{\ell} \simeq \ell_{\perp} / \ell$ and write $\G_{\ell} = \G
\cap N_{\ell}$. The boundary component corresponding to the cusp
$\ell$ is the $(p-1)$-torus $ \G_{\ell} \back  N_{\ell}$ with
universal cover $\ell^{\perp} / \ell$. We then define the $(p-n)$
cycle $B_{U}$ at the boundary by
\begin{equation}
     B_U = \{ w \in \ell^{\perp} / \ell : (U,w) = 0 \}.
\end{equation}

We write $C_U = \G_{U} \backslash B_U$  with $\G_U = N_U \cap
\G_{\ell}$ and note that in the Borel-Serre compactification this
cycle only depends on the equivalence class of the cusp $\ell$;
i.e., we have $C_U = C_{\g U}$ with $\g \in \G$ such that $\g \ell
= \ell_i$ for some $i$.

An orientation for $U$ gives one for $C_U$ in the following way:

Pick any null line $\ell'= \Q u'$ as above such that $\ell$ and
$\ell'$ span a hyperbolic plane whose orthogonal complement in $V$
we denote by $W$. Recall that the forward pointing vectors (see
(\ref{forward})) give an orientation for $\ell$ and $\ell'$
respectively. The orientation of $U$ induces one for $U'= U \cap
W$ by requiring that the orientation of $\ell$ followed by the one
of $U'$ gives the orientation for $U$. $B_U$ is isomorphic to the
orthogonal complement of $U$ in  $\ell_0 \perp W$, and we require
that the orientation of $\ell$ followed by the ones of $B_U =
T_Z(B_U)$, $U'$ and finally of $\ell'_0$ gives the orientation of
$V$.

A fixed orientation for $U$ defines  a sign character
$\epsilon(X)$ for $X$  a rational $n$-frame with $span_{\Q}(X)
=U$, by setting $\epsilon(X) =1$ if $X$ defines the same
orientation on $U$ and  $\epsilon(X) =- 1$ otherwise. So $C_X =
\epsilon(X) C_U$.

\begin{remark}\label{orientation}
   When working with coordinates for $B$ adopted to $U$ (see Section 2) one
   obtains a different orientation for $C_U$ which differs from the given one
   by a factor of $(-1)^{(n-1)(p-n)}$.
\end{remark}

\vspace{.5cm}

The construction of a composite cycle in this situation is more
complicated:

Let $\beta \in Sym_n(\Q)$ positive semidefinite and of rank $n-1$.
We define
\begin{equation}
\Omega_{\beta}^s = \{ X \in V^n : \frac12(X,X)=\beta \quad
\text{and rank}(X) =n \},
\end{equation}
the 'singular' part of the hyperboloid $\Omega_{\beta}$. Since
$\beta$ is singular, the radical $R(X)$ of the span of $X \in
\Omega^s_{\beta}$ is nonzero, i.e., $R(X) = \ell =\ell_X$  for
some rational isotropic line $\ell$. For such a line, we define
\begin{equation}
\Omega_{\beta,\ell} = \{ X \in V(\Q) : (X,X)=\beta \quad
\text{and} \quad R(X)= \ell \} \subset \ell^{\perp}.
\end{equation}

We then have
\begin{equation}
\Omega_{\beta}^s = \coprod_{j=0}^{r} \coprod_{\g \in \G_j \back
\G} \Omega_{\beta,\g^{-1}\ell_j},
\end{equation}
where $\ell_0, \cdots,\ell_r$ are the cusp representatives of the
$\G$-orbits of rational isotropic lines.

We also write $\Omega_{\beta,j}$ for  $\Omega_{\beta,\ell_j}$ and
\begin{equation}
\mathcal{L}_{\beta,j} = \Omega_{\beta,j} \cap L^n + h.
\end{equation}

\begin{lemma}
Let $\ell$ be a rational isotropic line. There is a finite number
of rational $n$-dimensional subspaces $U_1,\cdots,U_a$ of
$\ell^{\perp}$ such that
\begin{equation}
\{U_1,\cdots,U_a\} = \{ \text{span}(X) : X \in \Omega_{\beta,\ell}
\cap L^n + h \}.
\end{equation}
\end{lemma}

\begin{proof}
Indeed, we consider the quadratic space $\ell^{\perp} / \ell$
which is positive definite. Then there are only finitely many
$\bar{X} \in (\Omega_{\beta,\ell} \cap L^n + h)/ \ell$ such that
$\tfrac12 (\bar{X},\bar{X}) = \beta$. Pulling back to
$\ell^{\perp}$ then gives the lemma.
\end{proof}

For each cusp $\ell_j$, find a collection $U_{ij}, i = 1, \dots,
a_j$, of $n$-dimensional subspaces of $\ell_j^{\perp}$ as in the
lemma. We will write
\begin{equation}
    \Omega_{\beta,i,j} = \{ X \in \Omega_{\beta,j} : span(X) = U_{ij} \}
        \qquad \text{and} \qquad
    \mathcal{L}_{\beta,i,j} =  \Omega_{\beta,i,j} \cap  L^n + h,
\end{equation}
so that
\begin{equation}
     \Omega_{\beta,j} = \coprod_{i=1}^{a_j} \Omega_{\beta,i,j}
        \qquad \text{and} \qquad
      \mathcal{L}_{\beta,j} =  \coprod_{i=1}^{a_j} \mathcal{L}_{\beta,i,j}.
\end{equation}

\begin{lemma}\label{group action}
The action of $\G_j$ on $\Omega_{\beta,j}$ induces a free action
of $ \G_{U_{ij}} \back \G_j$ on $\Omega_{\beta,i,j}$. Here $
\G_{U_{ij}} = N_{U_{ij}} \cap \G_j.$
\end{lemma}

\begin{proof}
We show that the action of $\G_j$ on  $\Omega_{\beta,j}$ carries
$\Omega_{\beta,i,j}$ into itself and that the induced action of
$\G_{U_{ij}}$ is trivial.

Indeed, an element $\g \in \G_j$ operates on an element $x \in
\ell_j + W = \ell_j^{\perp}$ by adding a multiple of $u_j$ to $x$.
Thus $\g$ leaves stable any subspace of $(\ell_j)^{\perp}$
containing $\ell_j$, whence $\g$ leaves
 $\Omega_{\beta,i,j}$ stable. Consequently, $\G_j$ leaves  $\Omega_{\beta,i,j}$ stable.
Also, $\G_{U_{ij}}$ acts trivially on $U_{ij}$ whence it acts
trivially on
 $\Omega_{\beta,i,j}$.

Finally, if $\g \in \G_j$ satisfies $\g X =X$, then, since $X$
spans $U_{ij}$, necessarily $\g |_{U_{ij}} = 1$.
\end{proof}

Let $\mathcal{C}_{\beta,i,j}$ be a set of coset representatives of
this action on $\mathcal{L}_{\beta,i,j}$, i.e.,
\begin{equation}
\mathcal{C}_{\beta,i,j} = \left(  \G_j /  \G_{U_{ij}} \right)
\back \left( \Omega_{\beta,i,j}  \cap L^n + h \right).
\end{equation}
We will see below that $\mathcal{C}_{\beta,i,j}$ is infinite. It
is clear that the collection of the $\mathcal{C}_{\beta,i,j}$
provides a set of representatives for $\G \back \Omega_{\beta}
\cap L^n +h$.

Pick $ a \in \Q^n $ in the radical of $\beta$. Then for all $ X\in
\Omega_{\beta}^s $,
\begin{equation}
X \cdot a \in \ell_X = \Q u_X,
\end{equation}
with $u_X \in \ell_X$ as in (\ref{forward}). We can take $a$
nonzero and primitive in $\Z^n$. With this condition, $X$
determines $a$ up to $\pm 1$, and we write $X \cdot a = \nu(X)
u_X$, where $\nu(X)$ is determined up to a sign. Following
\cite{KShin} we call $X$ \emph{reduced} if with such a choice of
$a$ we have
\begin{equation}
X \cdot a = \nu(X) u_X
\end{equation}
with $\nu(X) \in [0,1)$. Note that if $X$ is reduced so is $\g X$
with $\g \in \G$ and $\nu(X) = \nu(\g X)$. We write
$\Omega^{red}_{\beta}$ for the set of reduced elements in
$\Omega_{\beta}^s $.

\begin{lemma}\label{cosetreps}
 $\G$ acts with finitely many orbits on the reduced elements in
 $\G \back \Omega_{\beta}^s \cap L^n +h$, and $\mathcal{C}^{red}_{\beta,i,j}:= \mathcal{C}_{\beta,i,j} \cap \Omega^{red}_{\beta}$ forms a set of representatives.
\end{lemma}

\begin{proof}
It is enough to show that for each pair $i,j$,
$\mathcal{C}^{red}_{\beta,i,j}$ consists of only finitely many
elements. We write  $\mathcal{C}_{\beta,U_{ij},h}$ for $
\mathcal{C}_{\beta,i,j}$.

Choose $m \in SL_n(\Z)$ such that $me_0=a$, where $e_0 =
{^t(1,0,\cdots,0)}$. We put $\beta_0 = {^tm \beta m}$ and $k=hm$.
Then right multiplication by $m$ gives a bijection from
$\mathcal{C}_{\beta,U_{ij},h}$ to
$\mathcal{C}_{\beta_0,U_{ij},k}$.

As the vector $e_0$ is a primitive integral vector in the radical
of $\beta_0$, we have
\begin{equation}\label{betaformula}
\beta_0  =
\begin{pmatrix}
0 & 0 \\
0 & \beta_0'
\end{pmatrix}.
\end{equation}
where $\beta_0'$ is a positive definite $(n-1)$ by $(n-1)$ matrix.
We write $ U_{ij} = \ell_j \perp U'$ with $U'$ positive definite.
Picking an appropriate basis for $U'$ we can assume that the $n$ by
$n$ matrix $g(Y)$ for $Y \in \mathcal{C}_{\beta_0,U_{ij},k}$ is of
the form
\begin{equation}\label{Yformula}
g(Y) =
\begin{pmatrix}
y_0 \\
0   \quad   Y_1'
\end{pmatrix}
=
\begin{pmatrix}
y_{01} & y_0' \\
0 &      Y_1'
\end{pmatrix},
\end{equation}
where $y_0 = (y_{01},y'_0)$ is a row vector of size $n$ and $Y_1'$
is an invertible $(n-1)$ by $(n-1)$ matrix. Similarly, the
congruence condition $k$ is of the form
 \begin{equation}
g(k) =
\begin{pmatrix}
k_{01} & k_0' \\
0 &      k_1'
\end{pmatrix}.
\end{equation}
Also, since $y_1 \equiv k_1 \mod{\Z u_0}$ we have $y_{01} \equiv
k_{01} \mod{\Z}$. Since $\beta'_0$ is positive definite, there are
only finitely many $Y'_1$ which represent $\beta'_0$. We have
\begin{equation}
g(n(u')Y) =
\begin{pmatrix}
y_{01} & y_0' +2(u',Y'_1) \\
0 &      Y_1'
\end{pmatrix}
\end{equation}
for $n(u') \in \G_j$. Hence there are only finitely many $y_0'$,
but $y_{01}$ runs through the set $\{k_{01} +n : n\in \Z, n \neq
-k_{01}\}$. Assuming $k_{01} \in [0,1)$ we observe $\nu(Y) =
k_{01}$ for $Y$ reduced. This proves the assertion.

\end{proof}

The proof of the lemma shows that the representatives of
$\mathcal{C}_{\beta,i,j}$ come in natural $\Z$-classes: If $ X \in
\Omega_{\beta,\ell}$ and $X \cdot a \in \ell_X$ as above, then
$\{\tilde{X} = X + u_X {^ta'}: k \in \Z\}$ with $a' \in \Z^n$ and
$^ta'a =1$ defines the $\Z$-class. From this we see that each
class contains exactly two reduced frames. If $X$ is reduced with
respect to $a$, then $\tilde{X} = X - u_X {^ta'}$ with $a' \in
\Z^n$ and $^ta'a =1$ is reduced with respect to $-a$. Moreover,
$\nu(\tilde{X}) = 1 -\nu(X)$.

Recall that the first periodic Bernoulli polynomial is defined by
\begin{equation}
\mathbf{B}_1(\alpha) =
\begin{cases}
\alpha - \tfrac12  \qquad \text{if} \quad  \alpha \in (0,1) \\
0 \qquad  \qquad \text{if} \quad \alpha = 0.
\end{cases}
\end{equation}

We readily check
\begin{equation}
\mathbf{B}_1(\nu(X))\epsilon(X) =
\mathbf{B}_1(\nu(\tilde{X}))\epsilon(\tilde{X}).
\end{equation}

We are finally ready to define the singular weighted composite
cycle $C^s_{\beta}$ by
\begin{equation}
C^s_{\beta}= \sum_{X \in \G \back \Omega^{red}_{\beta} \cap L^n
+h} \frac12 \mathbf{B}_1(\nu(X)) C_X.
\end{equation}

\begin{remark}

We could also define for a complex parameter $s$ the cycle
\begin{equation}
C_{\beta,s} =  \frac12 \sum_{X \in \G \back \Omega_{\beta} \cap
L^n +h} |\nu(X)|^{-s} C_X.
\end{equation}
Then the arguments of the previous lemma show that $C_{\beta,s}$
converges for $Re(s) > 1$ and has a meromorphic continuation the
whole complex plane. For the value at $s=0$ we have
\begin{equation}
C_{\beta,0} = C^s_{\beta}.
\end{equation}
(see also the proof of Prop. \ref{Dirichlet}.)

\end{remark}

\section{A Cohomology Class for the Weil Representation}

Recall that the metaplectic cover $G'= Mp(n,\R)$ of the symplectic
group $Sp(n,\R)$ is a central group extension
\begin{equation}
1 \longrightarrow \C^1 \longrightarrow Mp(n,\R) \longrightarrow
Sp(n,\R) \longrightarrow 1
\end{equation}
of $Sp(n,\R)$. Here $\C^1 = \{z \in \C : |z| =1 \}$. We fix a
splitting $Mp(n,\R) = Sp(n,\R) \times \C^1$ and denote by $K'
\subset Mp(n,\R)$ the inverse image of the standard maximal
subgroup
\begin{equation}
\left\{
\begin{pmatrix}
a&b \\-b&a
\end{pmatrix}
\; : \; a +ib \in U(n) \right\}
\end{equation}
of $Sp(n,\R)$. Then $K'$ admits a character $\det^{1/2}$; i.e.,
its square descends to the determinant character of $U(n)$.

$G \times G'$ acts on the Schwartz space $S(V(\R)^n)$ via (the
restriction of) the Weil representation $\omega =\omega_{V(\R)}$
associated to the additive character $t \longmapsto e(t) :=
\exp(2\pi it)$, see for example \cite{Weil}. Recall that the
action of $G'$ on $\psi \in S(V(\R)^n)$ is characterized by the
formulae
\begin{equation}
\omega \left( \left( \begin{smallmatrix} a&0\\0&^ta^{-1}
    \end{smallmatrix} \right) \right) \psi(X) =
 (\det a)^{m/2}
  \psi(Xa)
 \end{equation}
 for $a \in GL^+_n(\R)$;
\begin{equation}
\omega \left( \left( \begin{smallmatrix} 1&b\\0&1
    \end{smallmatrix} \right) \right) \psi(X) =
e^{\pi i tr(b(X,X))} \psi(X)
\end{equation}
for $b \in Sym_n(\R)$;
\begin{equation}
\omega \left( \left( \begin{smallmatrix} 0&1\\-1&0
    \end{smallmatrix} \right) \right) \psi(X) = \g  \hat{\psi}(X),
\end{equation}
where $\hat{\psi}$ is the Fourier transform of $\psi$ and $\g$ an
eighth root of unity.

The central $\C^1$ acts by
\begin{equation}
\omega((1,t))\psi = \begin{cases} t\psi \quad &\text{if $m$ is odd} \\
 \psi \quad &\text{if $m$ is even} \end{cases}
\end{equation}
for all $t \in \C^1$.

The group $G$ acts on $S(V(\R)^n)$ via
\begin{equation}
\omega(g) \psi(X) = \psi(g^{-1}X),
\end{equation}
which commutes with the action $G'$.

For $Z \in B$, we define the corresponding Gaussian by
\begin{equation}
\varphi_0(X,Z) = \exp(-\pi tr(X,X)_Z)
\end{equation}
and put $\varphi_0(X) = \varphi_0(X,Z_0)$. Note that
$\varphi_0(X,Z)$ is $G$-invariant; i.e.,
\begin{equation}
\varphi_0(gX,gZ) = \varphi_0(X,Z).
\end{equation}

The space of differential $n$-forms on $B$ is
\begin{equation}
\mathcal{A} ^n(B) \simeq \left[ C^{\infty}(G) \otimes \bigwedge
^n(\mathfrak{p}^{\ast}) \right]^K,
\end{equation}
where the isomorphism is given by evaluating at $Z_0$.

The main result of \cite{KMI} (cf. also \cite{KM90}), specialized
to our situation, is the construction of a certain differential
$n$-form of $B$ with values in the Schwartz space $S(V(\R)^n)$.

\begin{theorem}[\cite{KMI}] \label{KM-results}

For each $n$ with $0 \leq n \leq p$, there is a nonzero Schwartz
form
\begin{equation}
\varphi_n \in \left[ S(V(\R)^n) \otimes \mathcal{A}^{n}(B) \right]
^G \simeq
 \left[ S(V(\R))^n \otimes \bigwedge ^n(\mathfrak{p}^{\ast}) \right]^K,
\end{equation}
such that
\begin{itemize}
\item[(i)]
\begin{equation}
 d \varphi_n = 0; \notag
\end{equation}
i.e., for each $X \in V(\R)^n$, $\varphi_n(X)$ is a closed
$n$-form on $B$ which is $G_X$-invariant:
\begin{equation}
g^{\ast}\varphi_n (X) = \varphi_n (X)\notag
\end{equation}
for $g \in G_X$, the stabilizer of $X$ in $G$.
\item[(ii)]
The forms are compatible with the wedge product:
\begin{equation}
\varphi_{n_1} \wedge \varphi_{n_2} = \varphi_{n_1+n_2},\notag
\end{equation}
where $\varphi_n =0$ for $n > p$.
\item[(iii)]

Assume $U = U(X)$ for a linear independent, positive definite
$n$-frame $X$ in $V(\R)$. Then a Poincar\'e dual of $ C_U = \G_U
\back B_U $ is given by
\begin{equation}
 \left[ e^{\pi (X,X)} \sum_{\g \in \G_U \back \G}
\g^{\ast}\varphi_n ( X) \right].\notag
\end{equation}
\end{itemize}
\end{theorem}

In \cite{KMI,KMII} Poincar\'e dual form means the following: Let
$C \subset M = \G \back B$ be a cycle of dimension $n$. The $\eta$
is a Poincar\'e dual form of $C$ if
\begin{equation}\label{PDForm}
\int_C \omega = \int_M  \omega  \wedge \eta
\end{equation}
holds for all \emph{compactly supported} (or rapidly decreasing)
closed $n$-forms $\omega$.

\medskip

We now give some explicit formulae for the forms $\varphi_n \in
 \left[ S(V(\R))^n \otimes \bigwedge ^n(\mathfrak{p}^{\ast}) \right]^K$.

Via the basis $\{e_1,\cdots,e_p\}$ for $Z_0^{\perp}$ we identify
$\mathfrak{p}$ with $\R^{p}$. Then $\omega_i$ becomes the
functional on  $\mathfrak{p}$ which picks out the $i$-th
coordinate. For $X = (x_1,...,x_n) \in V(\R)^n \simeq M_{m,n}(\R)$
(w.r.t. the basis $\{e_1,\dots,e_{p+1}\}$), $m=p+1$, and for $1
\leq s \leq n$, we then define the $1$-form
\begin{equation}
\omega(s,X) = \sum_{i=1}^{p} x_{is}\omega_i.
\end{equation}
Note that $\omega(s,X)$ only depends on the $s$-th column vector
$x_s$ of $X$: $\omega(s,X) = \omega(s,x_s)$. We  set
\begin{align}
2^{-n/2} \varphi_n(X) &= \left(\bigwedge_{s=1}^n \omega(s,X) \right) \; \cdot \; \varphi_0(X) \\
             &= \varphi_1(x_1)\wedge \cdots \wedge \varphi_1(x_n)
\end{align}
with $\varphi_0(X) = \exp( -\pi tr(X,X)_0 )$, as before.

Note that this differs from the corresponding quantity in
\cite{KMI} by a factor of $2^{n/2}$.

We easily see
\begin{equation}
\varphi_n(X) = 2^{n/2} \sum_{1 \leq j_1 < \cdots < j_n \leq p}
P_{j_1, \cdots ,j_n}(X)  \; \exp( -\pi tr(X,X)_0 ) \;\otimes
\omega_{j_1} \wedge \cdots \wedge \omega_{j_n},
\end{equation}
where $ P_{j_1 \cdots j_n}(X)$ is the determinant of the $n$ by
$n$ matrix obtained from $X$ by removing all rows except the $j_1,
\cdots, j_n$. Occasionally we will write $\hat{X}$ for this
matrix, suppressing the coordinates. We write $\varphi_{j_1,
\cdots ,j_n}(X)  = P_{j_1 ,\cdots ,j_n}(X)  \; \exp( -\pi
tr(X,X)_0 )$.

Then it is easy to see that we have
\begin{equation}
\varphi_n(X)(W) = 2 ^{n/2} \det(X,W) \; \exp( -\pi tr(X,X)_0 )
\end{equation}
for $W \in T_{Z_0}(B)^n \simeq \mathfrak{p}^n \simeq
{(Z_0^{\perp})}^n$. Lemma \ref{invariantforms} gives

\begin{corollary}\label{FORM.}
\[
       \sigma^{\ast} \varphi_n (X) = 2^{n/2}
       \sum_{1 \leq j_1 < \cdots < j_n \leq p}
        P_{j_1,  \cdots ,j_n}(X) e^{-\pi (X,X)_0} \otimes
        \nu_{j_1} \wedge \cdots \wedge \nu_{j_n}
\]
\end{corollary}

We write $\varphi_n(X,Z)$ for the corresponding $n$-form on $B$;
for $g \in G_0(\R)$, we have per construction
\begin{equation}
\varphi_n(gX,gZ) = \varphi_n(X,Z),
\end{equation}
which also implies Th. \ref{KM-results} (i).

\vspace{.5cm}

Fundamental for the relationship to modular forms is

\begin{theorem}[\cite{KMI,KMII}] \hfill

$\varphi_n$ is an eigenvector of the maximal compact $K' \subset
Mp(n,\R)$ under the action of the Weil representation. We have
\[
\omega(k') \, \varphi_n \; = \; \det(k')^{m/2} \varphi_n
\]
for $k' \in K'$.
\end{theorem}

We denote by $\mathcal{L}_m$ the $G'$-homogeneous line bundle over
$G'/K'$ to the character $\det^{-m/2}$ of $K'$. Then the previous
theorem can reformulated as
\begin{align}
\varphi_n &\in \left[ \mathcal{L}_m \otimes S(V(\R)^n) \otimes \mathcal{A}^{n}(B) \right] ^{G \times G'} \\
          &\simeq \left[ \C_{\chi_m} \otimes  S(V(\R)^n) \otimes \
            \bigwedge ^n(\mathfrak{p}^{\ast}) \right]^{K\times K'},
\end{align}
where $\C_{\chi_m}$ is the one-dimensional module on which $K'$
acts via the character $\det^{-m/2}$.

For $\tau = u+ iv \in \h_n = \{ \tau \in Sym_n(\C) : Im(\tau) >0
\} \simeq G'/K'$, the Siegel space of genus $n$, we define in the
usual way
\begin{equation}
\varphi_n(\tau,X,Z) = \det(v)^{-m/4} \omega(g'_{\tau})
\varphi_n(X,Z).
\end{equation}

Here $g'_{\tau} \in Sp_n(\R)$ is a standard element carrying the
base point $iI_n \in \h_n$ to $\tau$; i.e.,
\begin{equation}
g'_{\tau} =
\begin{pmatrix}
v^{\frac12} & v^{-\frac12} u \\
0 & v^{-\frac12}
\end{pmatrix}
=
\begin{pmatrix}
1 &  u \\
0 & 1
\end{pmatrix}
\begin{pmatrix}
v^{\frac12} & 0 \\
0 & v^{-\frac12}
\end{pmatrix}.
\end{equation}

This is well defined, and we obtain
\begin{proposition}\label{phitauformula}
\[
 \varphi_n(\tau,X,Z)(W) = 2^{n/2}
\det (v)^{1/2} \det (X,W) e^{\pi i tr (X,X)_{\tau,Z}}
\]
for $W \in \left(T_Z(B)\right)^n \simeq \left(Z^{\perp}\right)^n$
and with $ (X,X)_{\tau,Z} = u(X,X) + iv(X,X)_Z $.
\end{proposition}

For a congruence condition $h \in (L^{\#})^n$, we define the theta
series $\theta(\tau)$ with  values in the differential $n$-forms
of $B$ by
\begin{equation}
\theta(\tau,Z) =  \sum_{X \in h + L^n} \varphi_n(\tau,X,Z).
\end{equation}

By the standard machinery of the theta correspondence (Poisson
summation formula) we get

\begin{theorem}[\cite{KMI,KMII}] \hfill

$\theta(\tau,Z)$ is a non-holomorphic Siegel modular form of
weight $m/2$ with values in the $\G$-invariant differential forms
of $B$ for some suitable congruence subgroup of $Sp(n,\Z)$.

\end{theorem}

\vspace{.5cm}

In \cite{KM90} it was shown that $\bar{\partial} \varphi_n$ (with
respect to the symplectic variable $\tau \in \h$) is \emph{exact}
in the orthogonal variable $Z \in B$; i.e., there exists
\begin{equation}
\psi_{n-1} \in [\mathcal{L}_m \otimes S(V(\R))^n \otimes
\mathcal{A}^{n-1}(B) \otimes \mathcal{A}^{0,1}(\h_n)]^{G \times
G'}
\end{equation}
such that
\begin{equation}
\bar{\partial} \varphi_n = d \psi_{n-1}.
\end{equation}

Defining the analogous theta series
\begin{equation}
 \Xi(\tau,Z) =  \theta_{\psi}(\tau,Z) =  \sum_{X \in h + L^n} \psi(\tau,X,Z)
\end{equation}
we  obtain
\begin{equation}
\bar{\partial} \; \theta(\tau,Z) = d \; \Xi(\tau,Z).
\end{equation}

\vspace{.5cm}

We now give a concrete formula for $\psi_{n-1}$. Consider the
double complex
\begin{equation}
 [\mathcal{L}_m \otimes S(V(\R))^n \otimes \mathcal{A}^{i}(B) \otimes \mathcal{A}^{0,j}(\h_n)]^{G \times G'}
\end{equation}
with maps $d, \bar{\partial}$. The Lie-algebra version of this
complex is the following. Let $\mathfrak{g}' = \mathfrak{k}' +
\mathfrak{p}'$ be the complexified Cartan decomposition of
$\mathfrak{sp}_n$. We can identify $ \mathfrak{p}'$ with the
complex tangent space of $\h_n \simeq G'/K'$ at the base point
$i1_n$, and the Harish-Chandra decomposition  $\mathfrak{p}'
=\mathfrak{p}_+ \oplus \mathfrak{p}_-$ gives the splitting of
$\mathfrak{p}'$ into the holomorphic and antiholomorphic tangent
spaces. We let $\nu_{jk}, 1 \leq j \leq k \leq n$ be dual to the
standard basis of $\mathfrak{p}_- \subset Sym_n(\C)$. Evaluation
at the base points gives an isomorphism of the above complex with
\begin{equation}
  C^{i,j} =  [ \C_{\chi_m} \otimes S(V(\R))^n \otimes \bigwedge^i {\mathfrak{p}}^{\ast} \otimes
\bigwedge^j \mathfrak{p}^{\ast}_-]^{K \times K'}.
\end{equation}
 We define $d, \bar{\partial}$ on $ C^{i,j}$ via transport of structure, for explicit formulae see \cite{KM90}. Note that
$\varphi_n \in C^{n,0}$ and $\psi_{n-1} \in C^{n-1,1}$.

We put (in coordinates for $\{e_i\}$)
\begin{equation}
A_{jk} (X) = (-1)^{k-1} x_{m,j} e^{-\frac12 \pi (x_k,x_k)_0}
\varphi_1(x_1) \wedge \cdots \wedge \widehat{\varphi_1(x_k)}
\wedge \cdots \wedge \varphi_1(x_n),
\end{equation}
where $\widehat{}$ over a term denotes that this term is omitted
in the product. We have
\begin{equation}
A_{jk} = (-1)^{k-1} \sum_{1 \leq \alpha_1 < \cdots \alpha_{n-1}
\leq p}  x_{mj} \, P^{(k)}_{\alpha_1, \cdots ,\alpha_{n-1}}(X)
\varphi_0(X) \, \otimes \omega_{\alpha_1} \wedge \cdots \wedge
\omega_{\alpha_{n-1}}
\end{equation}
Here $ P^{(k)}_{\alpha_1, \cdots ,\alpha_{n-1}}(X)$ is the
following polynomial. Let $X^{(k)}$ denote the $m$ by $(n-1)$
submatrix of $X$ obtained by deleting the $k$-th column. Then
 $ P^{(k)}_{\alpha_1, \cdots ,\alpha_{n-1}}(X) =  P^{(k)}_{\alpha_1, \cdots ,\alpha_{n-1}}(X^{(k)})$ is the minor obtained from $X^{(k)}$ using the rows ${\alpha_1, \cdots ,\alpha_{n-1}}$.

We now define
\begin{equation}
\psi_{n-1} = i 2^{n/2} \left[ \sum_{1 \leq j \leq n} A_{jj}  \otimes
\nu_{jj} \nu_{j_1} \wedge \cdots \wedge \nu_{j_n}; + \; \frac34
\sum_{1 \leq j < k \leq n} (A_{jk} + A_{kj}) \otimes \nu_{jk}
\right].
\end{equation}

Then
\begin{theorem}[\cite{KM90}]
\[
\bar{\partial} \, \varphi_n = d \, \psi_{n-1}.
\]
\end{theorem}

We write $\psi_{jj;\alpha_1,\cdots,\alpha_n}$ and
$\psi_{jk;\alpha_1 ,\cdots  ,\alpha_n}$ for the coefficient of
$\omega_{\alpha_1} \wedge \cdots \wedge \omega_{\alpha_{n-1}}$ in
$A_{jj}$ and $A_{jk} + A_{kj}$ respectively.

\section{The Growth of $\theta(\tau,Z)$ and $\Xi(\tau,Z)$}

In this section we prove that $\theta(\tau,Z)$ extends to the
Borel-Serre boundary and that $\Xi(\tau,Z)$ is rapidly decreasing
on $\G \back B$.

Since $\varphi_n(X,g^{-1}Z) = \varphi_n(gX,Z)$, it suffices to
prove the required estimates on the fixed Siegel set
$\mathfrak{S}'$. By some standard arguments we can also assume
that the lattice $L$ is of the form
\begin{equation}
L = L \cap \ell_0 \; + \; L \cap W \; + \; L \cap \ell_0'.
\label{decomposition}
\end{equation}

\vspace{.5cm}

We first consider an arbitrary $n$-form $\varphi\in \left[
S(V(\R)^n \otimes \bigwedge ^n \mathfrak{p}^{\ast} \right]^K$ in
the polynomial Fock space, that is, the space of Schwartz
functions of the form $p(X) \varphi_0(X)$ with $p$ a polynomial
function on $V(\R)^n$. (In the Fock model of the Weil
representation these become polynomials on $\C^{nm}$.)

We extend our basis $\omega_1, \cdots ,\omega_n$ of $\mathfrak{p}$
to a frame field $V_1(Z), \dots,V_p(Z)$. We then have
\begin{multline}
\theta_{\varphi}(\tau,Z(t,b)) \left( V_{i_1}(Z(t,b)), \cdots, V_{i_n}(Z(t,b)) \right) \\
= \sum_{X \in L^n +h} \varphi(\tau,Z_0,a(t)^{-1}n(b)^{-1}X)
(w_{i_1},\cdots,w_{i_n}).
\end{multline}

Thus the problem of estimating a form of the above type on
$\mathfrak{S}'$ reduces to estimating an expression of the
following type
\begin{equation}
\theta(t,b,R) = \sum_{X \in h + L^n} p\left( a(t)^{-1} n(b)^{-1} X
\right) \exp \left( -\pi R ( a(t)^{-1} n(b)^{-1} X ) \right),
\end{equation}
where $p(X)$ is a homogeneous polynomial function on $V^n$ and $R$
is a complex-valued quadratic function on $V^n$ with positive
definite real part. We now make some elementary observations
concerning the growth of such expressions in $t$. We define
$\theta^{\ast}(t,b,R)$ by
\begin{equation}
\theta^{\ast}(t,b,R) = \sum_{X \in h +L^n} \left| p\left(
a(t)^{-1} n(b)^{-1} X \right) \exp \left(- \pi R( a(t)^{-1}
n(b)^{-1} X ) \right) \right|.
\end{equation}
Via $V(\R)^n \simeq M_{m,n}(\R)$ we think of $p$ as a polynomial
in some coordinates of $V(\R)^n$.

From now on we use coordinates $y_{ij}$ with respect to a Witt
basis, see Section 2. Writing $X = \left( \begin{smallmatrix} y_1
\\ Y' \\ y_m \end{smallmatrix} \right)$ we have
\begin{equation}
a(t)^{-1} n(b)^{-1} X =
\begin{pmatrix}
t^{-1}\left( y_1 - 2(Y',b) + (b,b)y_m \right)\\
Y' - b \cdot y_m \\
ty_m
\end{pmatrix}.
\end{equation}
As a warm up we note
\begin{lemma} \label{warm up}
Suppose $p(y_{ij})$ is in the ideal in $\C[y_{ij}]$ generated by
$\{ y_{mj} : 1 \leq j \leq n \}$, the ideal of polynomial
functions which vanish on $(\ell_0^{\perp})^n$. Then
$\theta^{\ast}(t,b,R)$ is exponentially decreasing on
$\mathfrak{S}'$.
\end{lemma}

\begin{proof}
We may replace $R$ by $c \sum_{i,j} y_{ij}^2$ for a suitable
constant $c$ (since we are taking absolute values of the terms in
the sum). Under the hypothesis of the lemma the only terms that
appear in the sum have $y_{mj}(X) \neq 0$ for some $j$. But these
terms appear in the exponential multiplied by $t^2$, and the lemma
follows.
\end{proof}

\begin{remark}
The forms $\psi_{n-1}$ and $\varphi_n$ are \emph{not} of this
form.
\end{remark}

\vspace{.5cm}

We have an isomorphism
\begin{equation}
S(V(\R)^n) \longrightarrow S((\ell_0)^n) \otimes S(W(\R)^n)
\otimes S((\ell_0')^n)
\end{equation}
given by the partial Fourier transform operator
\begin{equation}
\mathcal{F}_0(\varphi_1 \otimes \varphi_2 \otimes \varphi_3) =
\widehat{\varphi_1} \otimes
 \varphi_2 \otimes \varphi_3.
\end{equation}
Here $  \widehat{\varphi_1}$ is the usual Fourier transform on
$(\ell_0)^n$. The right hand side is sometimes referred to as the
mixed model of the Weil representation.

We will need some formulae relating the action of $\omega$ and
$\mathcal{F}_0$ on $S(V(\R)^n)$. Identifying $(\ell_0)^n$ with
$M\Z^n$ for some $M \in \Q$ we denote the Fourier transform
variable (dual to $y_1$) by $\xi \in \R^n$.

\begin{lemma}\label{4.3}
Let $(\xi,w,y_m) \in (\R)^n \times W(\R)^n \times
(\ell_0'(\R))^n$.
\begin{itemize}
\item[(i)] For $n(b) \in N(\R)$ with $b \in W$,
      \begin{equation*}
      \mathcal{F}_0(n(b)\varphi)(\xi,w,y_m) =
              e \left( \xi^{t}(-2(b,w) + (b,b)  y_m)\right)
                       \mathcal{F}_0\varphi(\xi,w - by_m,y_m);
      \end{equation*}
\item[(ii)] For $a(t) \in A(\R)$,
      \begin{equation*}
       \mathcal{F}_0(a(t)\varphi) (\xi,w,y_m) = t^n \mathcal{F}_0\varphi(t\xi,w,ty_m);
      \end{equation*}
\item[(iii)] For $a'(v) = \left( \begin{smallmatrix} v&0\\0&^tv^{-1} \end{smallmatrix}  \right)
        \in Sp(n,\R)$ with $v \in GL_n(\R)$,
      \[
       \mathcal{F}_0(a'(v)\varphi)(\xi,w,y_m)  = (\det v)^{\frac{m}2 -1}
                                                 \mathcal{F}_0\varphi(\xi ^{t}v^{-1},wv,y_mv)
      \]
\item[(iv)] For $n'(u) = \left( \begin{smallmatrix} 1&u\\0&1 \end{smallmatrix}  \right)\in Sp(n,\R)$ with $u \in Sym_n(\R)$,
      \[
       \mathcal{F}_0(n'(u)\varphi)(\xi,w,y_m) = e \left(tr(u\frac{(w,w)}{2}) \right) \mathcal{F}_0\varphi(\xi+ \frac12 y_mu ,w         ,y_m).
       \]
\end{itemize}
\end{lemma}

\begin{proof}
This is an easy exercise which we omit.
\end{proof}

We introduce the following notation
\begin{align}
a^{\times}(v)\varphi(\xi,w,y_m) &= \varphi(\xi ^{t}v^{-1},wv,y_mv), \\
\phi(b,\xi,w,y_m) &=  e \left( \xi^{t}(-2(b,w) + (b,b)y_m)\right).
\end{align}
Note $|\phi(b,\xi,w,y_m)| = 1$.

Let $I \subset S(V^n)$ be the ideal of Schwartz functions in the
polynomial Fock space that vanish on the linear subspace $W^n$ of
$V^n$. Note that $W^n$ is defined by the equations
\begin{equation}
y_{1j} = 0 \qquad and \qquad y_{mj} = 0
\end{equation}
for $j= 1,\cdots,n$, i.e., $I = < y_{1j},y_{mj} >$. We observe
that if $\mathcal{F}_0 \varphi \in I$ then also $
\mathcal{F}_0(a'(v)\varphi)$ and $
 \mathcal{F}_0(n'(u)\varphi) \in I$.

\begin{lemma}\label{4.2}
Suppose $\mathcal{F}_0 \varphi$ is in the  ideal $I$. Then
$\theta(t,b)$ is exponentially decreasing on $\mathfrak{S}'$.
\end{lemma}

\begin{proof}
We may write
\begin{equation}
\theta(t,b) = < \Theta_{h +L^n}, n(b)a(t)\varphi >.
\end{equation}
Here $\Theta_{h+L^n}$ is the sum of Dirac deltas (placed at the
points of $h+L^n$) and $<\;,\;>$ denotes the Kronecker pairing.

We write $h = h_1 + h'$ with $h_1 \in (\ell_0)^n$. Then there is a
constant $C$ such that
\begin{equation}
\mathcal{F}_0 \Theta_{h+L^n} = C e (\xi ^{t}h_1) \Theta_{h'+L^n}
\end{equation}
by Poisson summation. Using the formulas from the previous lemma
we obtain
\begin{equation}
\mathcal{F}_0( (n(b)a(t)) \varphi)(\xi,w,y_m) = \phi(b,\xi,w,y_m)
\mathcal{F}_0\varphi(t\xi,w,ty_m).
\end{equation}
Hence
\begin{equation}
\theta(t,b) = C t^n \sum_{ \substack{ \xi \in M^{-1}\Z^n \\
(w,y_m) \in W \times (\ell_0')^n +h'}}  \phi(b,\xi,w,y_m) e(\xi \,
^th_1) \mathcal{F}_0\varphi(t\xi,w,ty_m).
\end{equation}
The lemma now follows from an argument analogous to that of Lemma
\ref{warm up}.
\end{proof}

Note however that one cannot conclude from the lemma that
$\theta^{\ast}(t,b)$ is rapidly decreasing.

\begin{lemma}\label{4.5}
Suppose $p(y_{ij})$ is divisible by $y_{1j}$ for some $j$ but no
higher power of $y_{1j}$. Then $\mathcal{F}_0 \varphi$ is in the
ideal $I$. Moreover, for every $v \in GL_n(\R)$ and $u \in
Sym_n(\R)$ the function $\mathcal{F}_0(n'(u)a'(v) \varphi)$ is in
the ideal $I$.
\end{lemma}

\begin{proof}
The first statement is clear for we may write
\begin{equation}
\varphi(X) = y_{1j} e^{-\pi y_{1j}^2} \psi(X).
\end{equation}
where $\psi(X)$ does not involve $y_{1j}$. Now taking
$\mathcal{F}_0$ does not change the function since $ y_{1j}
e^{-\pi y_{1j}^2}$ is its own Fourier transform (up to the
constant $-i$).

The second statement follows from the first one, the formulae
(iii) and (iv) of Lemma \ref{4.3} and the observation that if
$\mathcal{F} \varphi \in I$, then $a^{\times}(v) \mathcal{F} \varphi \in
I$ and \linebreak $e \left(tr(u\tfrac{(w,w)}{2}) \right)
\mathcal{F}_0\varphi(\xi+ \tfrac12{y_mu},w,y_m) \in I$.
\end{proof}

\begin{corollary}
If $\varphi$ satisfies the hypothesis of the lemma, then for every
$\tau \in \h_n$, $\theta(t,b,\tau)$ decays exponentially on
$\mathfrak{S}'$, where
\begin{equation}
\theta(t,b,\tau) = \sum_{X \in L^n +h} p(X v^{\frac12})
e^{-\pi(X,X)_{\tau,Z(t,b)}}.
\end{equation}
\end{corollary}

\vspace{.5cm}

We now check that the form $\psi_{n-1}$ satisfies the hypothesis
of the previous lemma. It is enough to do this for the individual
components $\psi_{jk;\alpha_1,\cdots,\alpha_{n-1}}$

\begin{lemma}\label{4.6} \hfill
\begin{itemize}
     \item[(i)]
        $\mathcal{F}_0   \psi_{jj;\alpha_1,\cdots,\alpha_{n-1}}       \in I \qquad \text{ for all} \qquad        1 \leq \alpha_1 < \cdots < \alpha_{n-1} \leq p$
     \item[(ii)]
        $\mathcal{F}_0   \psi_{jk;\alpha_1,\cdots,\alpha_{n-1}}     \in I \qquad \text{ for all} \qquad
        1 \leq \alpha_1 < \cdots < \alpha_{n-1} \leq p$
\end{itemize}
\end{lemma}

\begin{proof}
(i) follows immediately from Lemma \ref{4.5} and the explicit
formulae for $\psi_{n-1}$: We have
\begin{align}
\psi_{jj;\alpha_1,\cdots,\alpha_{n-1}}(X)
&= x_{mj} P^{(j)}_{\alpha_1,\cdots,\alpha_{n-1}}(X) \varphi_0(X) \\
&= (y_{1j} - y_{mj})  P^{(j)}_{\alpha_1,\cdots,\alpha_{n-1}}(X)
\varphi_0(X).
\end{align}
Now observe that $ P^{(j)}_{\alpha_1,\cdots,\alpha_{n-1}}(X)$ is a
polynomial which does not involve $y_{1j}$.

(ii) is more complicated. By an argument similar to the previous
one we find that $\mathcal{F}_0
\psi_{jk;\alpha_1,\cdots,\alpha_{n-1}}  \in I$ provided $\alpha_1
\neq 1$. However, for $\alpha_1 =1 $ it is no longer true. We
have, assuming $j <k$,
\begin{equation}
\psi_{jk; \alpha_1,\cdots,\alpha_{n-1}} =
           \left[ (-1)^{k-1} x_{mj}  P^{(k)}_{1,\cdots,\alpha_{n-1}}(X) + (-1)^{j-1} x_{mk}
                P^{(j)}_{1,\cdots,\alpha_{n-1}}(X) \right] \varphi_0(X).
\end{equation}
We expand $ P^{(k)}_{1,\cdots,\alpha_{n-1}}(X)$ by the first row
and obtain
\begin{equation}
 P^{(k)}_{1,\cdots,\alpha_{n-1}}(X) = (-1)^{j-1} x_{1j} P^{(j,k)}_{\alpha_2,\cdots,\alpha_{n-1}}(X) + R_k,
\end{equation}
where $R_k$ is a polynomial not involving $x_{1j}$. Similarly,
\begin{equation}
 P^{(j)}_{1,\cdots,\alpha_{n-1}}(X) = (-1)^{k-2} x_{1k} P^{(j,k)}_{\alpha_2,\cdots,\alpha_{n-1}}(X) + R_j.
\end{equation}

We obtain
\begin{equation}
\psi_{jk}(X) \equiv 2 (-1)^{j+k} (y_{1j}^2 - y_{1k}^2)
 P^{(j,k)}_{\alpha_2,\cdots,\alpha_{n-1}}(X) \varphi_0(X) \mod{\mathcal{F}_0^{-1} I}.
\end{equation}
Taking $\mathcal{F}_0$ now shows that the right hand side is in
$\mathcal{F}_0^{-1}I$. This boils down to the fact that the
Fourier transform of $x^2 e^{-\pi x^2}$ is $(
\frac{1}{2\pi}-x^2)e^{-\pi x^2}$.
\end{proof}

We conclude

\begin{theorem} \label{Xi-growth} \hfill

$\Xi(\tau,Z)$ is rapidly decreasing.
\end{theorem}


\vspace{.5cm}

We now determine the growth of $\theta(\tau,Z) =
\theta_{\varphi_n}(\tau,Z)$ on $\mathfrak{S}'$. Recall
\begin{equation}
\varphi_{\alpha_1,\cdots,\alpha_{n}}(X) = 2^{n/2}
P_{\alpha_1,\cdots,\alpha_{n}}(X) \varphi_0(X)
\end{equation}
and put $\theta_J(\tau,Z) = \sum_{X \in L^n +h}
\varphi_{\alpha_1,\cdots,\alpha_{n}}(X)$ with $J=
\{\alpha_1,\cdots,\alpha_{n}\}$.

We write
\begin{equation}
L^n +h = (\ell_0)^n + h_0) + (L^n \cap W^n + h_W) + (\ell_0')^n +
h'_0) \label{decomposition1}
\end{equation}
according to the decomposition (\ref{decomposition}).

The following lemma gives the growth of the components
$\theta_J(\tau,Z)$ of $\theta(\tau,Z)$.

\begin{lemma} \hfill

\begin{itemize}\label{estimate}
\item[(i)]
     $\varphi_{\alpha_1,\cdots,\alpha_{n}} \in \mathcal{F}_0^{-1} I$ if and only if $\alpha_1 =1$.
\item[(ii)]
     If $\alpha_1 =1$, then $\theta_J(\tau,Z)$ has exponential decay on $\mathfrak{S}'$.
\item[(iii)]
     If $\alpha_1 \neq 1$,  then $\theta_J(\tau,Z) =
     \begin{cases}
     O(t^n) \qquad \qquad \; \text{if} \quad h'_0 \in  (\ell_0')^n \\
     O(t^n e^{-C t^2}) \qquad \text{if} \quad h'_0 \notin (\ell_0')^n.
     \end{cases}$ \\ as $t \to \infty$.
\item[(iv)]
     $\theta_J^{\ast}(\tau,Z) = O(t^n)$ if  $h'_0 \in  (\ell_0')^n$.
\end{itemize}
\end{lemma}

\begin{proof}
For (i) develop $P_{\alpha_1,\cdots,\alpha_{n}}(X)$ after the
first row and proceed as in the proof of Lemma \ref{4.6}. (ii)
follows from (i) and Lemma \ref{4.2}. For (iii), we first write
\linebreak
$\varphi_{\alpha_1,\cdots,\alpha_{n}}(a(t)^{-1}n(b)^{-1}X) =
2^{n/2} P_{\alpha_1,\cdots,\alpha_{n}}(X) e^{-2\pi t^{-2} (\sum
y^2_{1k}) -\pi tr(X',X') - 2\pi t^2(\sum y^2_{mk}) }$ with
\linebreak $n(b)^{-1}X = (y_1,X',y_m)$. The assertion now follows
from $\sum_{k \in \Z+h} e^{-\pi (k/t)^2}= O(t)$ as $t \to \infty$,
which can be most easily seen by taking the Fourier transform, and
$\sum_{k \in \Z+h'} e^{-\pi (kt)^2}= O(e^{-Ct^2})$ if and only if
$h' \notin \Z$. This also implies $(iv)$ in the case of $\alpha_1
\neq 1$. If $\alpha_1 = 1$, then (iv) reduces to $\sum_{k \in
\Z+h} |\tfrac{k}{t}| e^{-\pi (k/t)^2}= O(t)$, which is an easy
calculus exercise.
\end{proof}

Note that the condition $ h'_0 \in   (\ell_0')^n$ certainly is
equivalent to $h \in  (\ell_0^{\perp})^n$. Following
(\cite{KShin}) we call the congruence condition $h \in (L^{\#})^n$
\emph{non-singular} if for all frames $X \in h+ L^n$ of rank $n$,
the radical $R(X)$ is empty. Otherwise we call $h$ singular.

Near the cusp given by $\ell_0$ we can change the upper-half space
coordinates $(t,b)$ to $(s,b)$ with $s =1/t$. Then the restriction
of a differential form on to the boundary component coming from
$\ell_0$ is given by setting $s=0$ (and corresponds to $t \to
\infty$).

\begin{theorem}\label{growth}(Theorem \ref{THA}) \hfill
\begin{itemize}
     \item[(i)]  $\theta(\tau)$ extends to the Borel-Serre boundary of $M$; i.e., defines a closed differential form on $\overline{M}$.

     \item[(ii)] If $h$ is non-singular or $n=p$, then $ \theta(\tau)$ is
                rapidly decreasing on $M$; hence
                $\theta(\tau)|_{\partial \overline{M}} = 0$.
     \item[(iii)] If $h$ is singular, then the restriction of $\theta_{\varphi}$
                 to the component of the Borel-Serre boundary $e_P$ coming from the                  parabolic $P$ is the restriction of the theta series to the
                 positive definite subspace $W$ of $V$. More precisely, under the
                 assumption (\ref{decomposition1}),
\[
 \theta_{\varphi}|_{e_P}(\tau,Z(b)) =
\begin{cases}
\underset{X \in  L^n \cap W^n + h_W}{\sum} \varphi(\tau,X,Z(b))  &\text{if} \quad h \in (\ell_0^{\perp})^n \\
\qquad 0  & \text{if} \quad h \notin (\ell_0^{\perp})^n.
\end{cases}
\]
Here
\[
\varphi(\tau,X,Z(b)) = \sum_{2 \leq \alpha_1 < \cdots < \alpha_{n}
\leq p} P_{\alpha_1 , \cdots , \alpha_{n}}(X) \exp(-\pi
tr(X,X)\tau) db_{\alpha_1} \wedge \cdots \wedge db_{\alpha_{n}}
\]
for $X \in W^n$.
\end{itemize}

\end{theorem}
\begin{proof}
Everything follows from Corollary \ref{FORM.},
Lemma \ref{invariantforms}, Lemma \ref{estimate} and
\begin{equation}
\lim_{t \to \infty } t^{-n} \theta_{  \alpha_1 , \cdots ,
\alpha_{n} }(\tau,Z)= P_{\alpha_1 , \cdots , \alpha_{n}}(X)
\exp(-\pi tr(X,X)\tau) db_{\alpha_1} \wedge \cdots \wedge
db_{\alpha_{n}},
\end{equation}
which is seen by taking the operator $\mathcal{F}_0$, Lemma
\ref{4.3} and Poisson summation.
\end{proof}

\begin{remark}
Theorem \ref{growth} (iii) also shows a nice functorial property
of the Weil representation. We have
\[
\left( \omega_{V(\R)}(g'(\tau) \theta_{\varphi}) \right)|_{e_P}
(Z(b)) =  \omega_{W(\R)} (g'(\tau)) \left(\theta_{\varphi}|_W
\right) (Z(b)),
\]
where $\omega_{W(\R)}$ is the Weil representation attached to the
positive definite space $W(\R)$ and $\theta_{\varphi}|_W $ is the
theta series restricted to $W$.
\end{remark}

\vspace{.5cm}

By Theorem \ref{growth}  we can now define for a closed
differential $(p-n)$-form $\eta$ on $\overline{M}$,
\begin{equation}
\Lambda(\eta)(\tau) = \int_{M} \eta(Z) \wedge \theta(\tau,Z).
\end{equation}
This extends the lift considered in \cite{KM90} to forms which do
not vanish at the boundary.

\begin{theorem}\label{holomorph}(Theorem \ref{THB}) \hfill

$\Lambda(\eta)(\tau)$ is a \underline{\emph{holomorphic}} Siegel
modular form of weight $m/2$.
\end{theorem}

\begin{proof}
This will now follow from Theorem \ref{Xi-growth} and the
following calculation, see \cite{KM90}:
\begin{multline}
    \bar{\partial} \Lambda(\eta)(\tau) =
    \bar{\partial}\int_{M} \eta(Z) \wedge \theta_{\varphi}(\tau,Z) =
    \int_{M} \eta(Z) \wedge \theta_{ \bar{\partial} \varphi}(\tau,Z)\\
    =\int_{M} \eta(Z) \wedge \theta_{d \psi }(\tau,Z) =
     \int_{M} d( \eta(Z) \wedge \theta_{ \psi }(\tau,Z)) = 0.
\end{multline}
The last equation is Stokes' Theorem. Here we need $\theta_{\psi}$
rapidly decreasing.
\end{proof}

\begin{remark}
In the analogous situation of locally symmetric spaces associated
to orthogonal groups of arbitrary signature $\theta_{\psi}$ is not
rapidly decreasing and the above argument breaks down. The theta
integral is non-holomorphic in general, see \cite{Funke}. In
\cite{KM90} it was assumed that $\eta$ was rapidly decreasing and
the above argument showed the holomorphicity of the theta
integral.
\end{remark}

\section{The Singular Fourier Coefficients}

In this section we compute the singular Fourier coefficients of
the theta integral $\Lambda(\eta)(\tau)$.

For the $\beta$-th Fourier coefficient, we have
\begin{equation}
a_{\beta}\left( \eta \right) = \int_{M} \eta \wedge \sum_{X \in
\Omega_{\beta} \cap (h + L^n)} \varphi(iv,Z,X) e^{-2\pi tr(\beta
v)}.
\end{equation}

First note that Prop. \ref{phitauformula} implies that only for
rank $(X) =n $ we have $\varphi(X) \neq 0$. Therefore $a_{\beta} =
0 $ unless $\beta = \tfrac12(X,X)$ is singular and positive
semidefinite with rank$(X) =n$ (or $\beta$ positive definite).

\begin{theorem}\label{semidef}(Theorem \ref{THD})
Assume that $\beta$ is positive semi-definite of rank $n-1$. Then
\[
a_{\beta}(\eta) = (-1)^n \int_{C_{\beta}} \eta.
\]
\end{theorem}

\begin{proof}
With the notation of Section 3 we have
\begin{align}
 e^{2\pi  tr(\beta v)} a_{\beta} &=
    \int_{\G \back B} \eta  \wedge \sum_{X \in \Omega_{\beta}^s \cap (h + L^n)}
                     \varphi(iv,Z,X)   \\
&=      \int_{\G \back B} \eta \wedge \sum_{j=0}^r \sum_{\g \in \G_j \back \G}           \sum_{X \in  \mathcal{L}_{\beta,j}} \g^{\ast} \varphi(iv,Z,X) \\
&=  \sum_{j=0}^r \int_{\G \back B} \eta \wedge \sum_{\g \in \G_j \back \G} \sum_{X \in    \mathcal{L}_{\beta,j}} \g^{\ast}\varphi(iv,Z,X) \\
&= \sum_{j=0}^r \sum_{i=1}^{a_j} \int_{\G \back B} \eta \wedge
\sum_{\g \in \G_j \back \G} \sum_{X \in \mathcal{L}_{\beta,i,j}}
\g^{\ast}\varphi(iv,Z,X).
\end{align}

\begin{proposition}\label{unfolding}
\[
     \int_{\G \back B} \eta \wedge \sum_{\g \in \G_j \back \G} \sum_{X \in    \mathcal       {L}_{\beta,i,j}} \g^{\ast}\varphi(iv,Z,X) =
      \int_{\G_j \back B} \eta \wedge \sum_{X \in    \mathcal{L}_{\beta,i,j}}
      \varphi(iv,Z,X).
\]
\end{proposition}

\begin{proof}
The considerations in Section 5 imply that it is enough to show
that \linebreak $\sum_{X \in \mathcal{L}_{\beta,i,j}} \varphi(iv,
a(t)^{-1} n(b)^{-1} X)$ is rapidly decreasing for $t \to \infty$
and $t \to 0$. Taking $m \in SL_n(\Z)$ as in the proof of Lemma
\ref{cosetreps} we find via Prop. \ref{phitauformula} that
\begin{equation}
\sum_{X \in    \mathcal{L}_{\beta,U_{ij},h}} \varphi(iv, a(t)^{-1}
n(b)^{-1} X) = \sum_{Y \in    \mathcal{L}_{\beta_0,U_{ij},k}}
\varphi(iv', a(t)^{-1} n(b)^{-1} Y)
\end{equation}
with $v'= {^t}m^{-1}vm^{-1}$ and $\beta_0$, $k=hm$ and $\beta$ as
in (\ref{betaformula}). But now Lemma \ref{phi-formula} and Lemma
\ref{estimate}(i) show that we have
$\varphi|_{\mathcal{L}_{\beta_0,U_{ij},k}}  \in \mathcal{F}_0^{-1}
I$ (in the notation of Section 5). As in (\ref{Yformula}), $Y \in
\mathcal{L}_{\beta_0,U_{ij},k}$ is of the form $ \left(
\begin{smallmatrix} y_0 \\ 0 \quad Y_1' \end{smallmatrix}\right)
$ with $y_0 \in (\ell_0)^n$ and finitely many possibilities for
$Y'_1$. So we can apply $\mathcal{F}_0$ to the sum over
$\mathcal{L}_{\beta_0,U_{i,j},k}$ and Lemma \ref{4.5} gives the
rapid decay as $t \to \infty$. The decay as $t \to 0$ is clear.
\end{proof}

\begin{remark}
In \cite{KShin} and \cite{KM90} unfolding was not attempted in the
above situation. This led to considerable complications. In
\cite{KShin}, the case $n=p$, Kudla introduced a wave packet
attached to the standard Eisenstein series for $O(p,1)$ to compute
the integral. The method employed in the following is
conceptually much simpler (even though the actual calculations
are quite similar). Moreover, it should be immediately available
in the more general situation of \cite{KM90} for not rapidly
decreasing $\eta$.
\end{remark}

We define a smooth differential $n$-form $\theta_{i,j}(Z)$ on
$\G_j \back B$ by
\begin{equation}
\theta_{i,j}(Z) = \sum_{X \in \mathcal{L}_{\beta,i,j}}
\varphi(Z,X).
\end{equation}
and put
\begin{equation}
\theta(\eta,\beta,U_{ij}) =
 \int_{\G_j \back B}  \eta \wedge \theta_{i,j}(Z).
\end{equation}
Hence
\begin{equation}\label{formula1}
 e^{2\pi  tr(\beta v)} a_{\beta} =
\sum_{j=0}^r \sum_{i=1}^{a_j} \theta(\eta,\beta,U_{ij}).
\end{equation}

We also define a function $\Phi(X,Z)$ via
\begin{equation}
\eta \wedge \varphi(Z,X) = \Phi(X,Z) d\mu,
\end{equation}
where $d\mu = t^{-p} dt\wedge db_1 \wedge \cdots \wedge db_{p-1}$
is the Riemannian volume form, and set
\begin{equation}
\Phi_{i,j}(Z) =  \sum_{X \in \Omega_{\beta,i,j}} \Phi(Z,X).
\end{equation}

Picking the standard fundamental domain for $\G_j \back B$ we
obtain
\begin{equation}\label{integral99}
\theta(\eta,\beta,U_{ij}) = (-1)^{p-1} \int_0^{\infty} \left(
\int_{\R^{p-1}/\Lambda_j}
  \Phi_{i,j}(Z(t,b))  db \right) t^{-p} dt.
\end{equation}

Recall (\ref{lattices}) that the lattice $\La_j \subset W \simeq
N_j \simeq \R^{p-1}$ is given by  $\La_j \simeq \G_j$. Here and
from now on $N_j = N_j(\R)$ and $W= W(\R)$. We denote the torus $
\Lambda_j \back W $ by $\mathbb{T}_j$.

We write $A_{i,j}(t)$ for the inner integral of
(\ref{integral99}). We have
\begin{align}
A_{i,j}(t) &= \int_{\G_j \back N_j} \Phi_{i,j}(Z(t,b))db \\
         &= \int_{\G_j \back N_j}  \sum_{X \in \mathcal{L}_{\beta,i,j}} \Phi(Z(t,b),X) db \\
         &=   \sum_{X \in \mathcal{L}_{\beta,i,j}} \int_{\G_j \back N_j} \label{covering} \Phi(Z(t,b),X) db.
\end{align}

We also get a splitting
\begin{equation}
W = U_{ij} \cap W + (U_{ij}^{\perp} \cap W). \label{splitting}
\end{equation}

Now note that the right hand side of (\ref{covering}) is
multiplicative under finite coverings! We pass to a subgroup
$\tilde{\G}_j \simeq \tilde{\Lambda}_j$ of finite index $\kappa$
in $\G_j$ given by
\begin{equation}
\tilde{\La}_j = (U_{ij} \cap \tilde{\La}_j) + (U_{ij}^{\perp} \cap
\tilde{\La}_j) =: \La' + \La''
\end{equation}
and obtain a degree $\kappa$ covering
\begin{equation}
     \tilde{\mathbb{T}}_j =\mathbb{T}' \times \mathbb{T}_{U_{ij}} \to \mathbb{T}_j
\end{equation}
with $\mathbb{T}'=  \La' \back U_{ij}$ and $\mathbb{T}_{U_{ij}} =
\La'' \back U_{ij}^{\perp}$. We obtain
\begin{equation}
           \kappa A_{i,j}(t) =
              \sum_{X \in \mathcal{L}_{\beta,i,j}}
              \int_{\mathbb{T}' \times \mathbb{T}_{U_{ij}}} \Phi(Z(t,b),X) db.
\end{equation}

We will use the horospherical coordinates adopted to $U_{ij}$. We
may choose a Witt basis $u,w_1,\dots,w_{p-1},u'$ with
$u,w_1,\dots,w_{n-1} \in U_{ij}$. Then the horospherical
coordinates $(t,b_1,\cdots,b_{p-1})$  are adopted to $U_{ij}$. The
decomposition
\begin{equation}
       Z(t,b) = Z(t,b',b'')
\end{equation}
corresponds to the splitting (\ref{splitting}). We write $\eta$ in
terms of these coordinates as
\begin{equation}\label{etaform}
          \eta(t,b) = f(t,b) db_n \wedge \cdots \wedge db_{p-1} + \eta'(t,b),
\end{equation}
where $\eta'(t,b)$ is in the ideal of forms on $\G_j \back B$
generated by $\{ dt,db_1,\cdots,db_{n-1} \}$.

For $\varphi_n(Z(t,b',b''),X)$, we have

\begin{lemma}\label{phi-formula}
Suppose $U := span(X) = span\{u_0,w_1,\cdots,w_{n-1}\}$. Then
\begin{equation*}
\varphi_n(Z(t,b',b''),X) = 2^{n/2} \frac12 \det g(X)
e^{-\pi(X,X)_{Z(t,b',0)}} \, t^{-n-1} dt \wedge db_1 \wedge \cdots
\wedge db_{n-1},
\end{equation*}
where $g(X)$ is the matrix expressing the basis $X$ for $U$ in
terms of the basis \linebreak $\{u_0,w_1,\cdots,w_{n-1}\}$ for
$U$; i.e.,
\[
(x_1,\cdots,x_n) =  (u_0,w_1, \cdots, w_{n-1}) g(X).
\]
\end{lemma}

\begin{proof}

 We have $\hat{X} = (x_{ij})$ with
\begin{equation}
x_{ij} =  (x_j,e_i) \qquad \text{for} \quad 1 \leq i \leq p \qquad
\text{and} \qquad x_{p+1,j} = - (x_j,e_{p+1}).
\end{equation}
Moreover $u_0 = \tfrac12 (e_1 + e_{p+1})$. But by assumption
$(x_j,u_0) =0$ for all $j$ and $(x_j,e_i) = 0$ for $i \geq n$. It
follows immediately  that the only vanishing Pl\"ucker coordinate
$P_{j_1,\cdots,j_n}(X)$, which is non-zero, is $P_{1, 2, \cdots,
n}$, and this has value $\tfrac12 \det g(X)$.

We next observe
\begin{equation}
\varphi_n (Z(t,b),X) = \varphi_n(Z_0, a(t)^{-1}n(b)^{-1}X).
\end{equation}
Now we have (using the previous formulas)
\begin{align}
\det g \left( a(t)^{-1}n(b)^{-1}X \right) &= t^{-1} \det g(X), \\
\exp(-\pi (X,X)_{Z(t,b)}) &= \exp (-\pi (X,X)_{Z(t,b',0)}), \\
\intertext{and} \nu_1 \wedge \cdots \wedge v_n &= t^{-n} dt \wedge
db_1 \wedge \cdots \wedge db_{n-1}.
\end{align}
The lemma now follows from Corollary \ref{FORM.}

\end{proof}

Writing $h(X,t,b') = \det(v)^{1/2} 2^{n/2} \frac12 \det g(X)
t^{-1} \varphi_0( iv,a(t)^{-1} n(b')^{-1}X)$, Lemma
\ref{phi-formula} and (\ref{etaform}) give
\begin{align}
     \eta \wedge    \varphi(iv,  Z(t,b',b'')) &=
     (-1)^{(p-n)n}   t^{-n} h(X,t,b') f(t,b',b'')
      dt \wedge db_1 \wedge \cdots \wedge db_{p-1} \\
\intertext{and}
     \Phi(X,Z(t,b',b'')) &=  (-1)^{(p-n)n}  t^{p-n} h(X,t,b') f(t,b',b'').
\end{align}

Thus the inner integral in (\ref{integral99}) is given by
\begin{equation}
        \kappa A_{i,j}(t) = {(-1)^{(p-n)n}} t^{p-n}
               \sum_{X \in \mathcal{L}_{\beta,i,j}}
               \int_{\mathbb{T}'} h(X,t,b')
               \left(\int_{\mathbb{T}_{U_{ij}}} f(t,b',b'')db''\right) db'.
\end{equation}
But the inner integral is equal to the period of the differential
form $\eta$ over the closed cycle $C_{U_{ij}}(t,b') \subset \G_j
\back B$ given by
\begin{equation}
     C_{U_{ij}}(t,b') = (-1)^{(n-1)(p-n)}n(b')a(t){\mathbb{T}_{U_{ij}}}.
\end{equation}
(For the sign, see Remark \ref{orientation}). Since $\eta$ is
closed, the period is independent of $b'$ and $t$, and we obtain
\begin{multline}
     \kappa A_{i,j}(t) =  2^{n/2}\det(v)^{\frac12} {(-1)^{p-n}}
     t^{p-n-1} \left( \int_{C_{U_{ij}}} \eta \right) \\
     \times \sum_{X \in \mathcal{L}_{\beta,i,j}}  \det g(X) \int_{\mathbb{T}'}      \varphi_0( iv, a(t)^{-1} n(b')^{-1}X) db'.
\end{multline}

We now unfold
\begin{equation}
    I = \sum_{X \in \mathcal{L}_{\beta,i,j}} \int_{\mathbb{T}'}
         \det g(X) \varphi_0(iv,a(t)^{-1}n(b')^{-1}X) db'.
\end{equation}
We observe $\mathbb{T}' \simeq \tilde{\mathbb{T}}_j/
\mathbb{T}_{U_{ij}} \to \mathbb{T}_j/ \mathbb{T}_{U_{ij}} \simeq
\tfrac{N/N_{U_{ij}}}{\G_j/\G_{U_{ij}}}$ is a covering of degree
$\kappa$.

We let $\mathcal{D}''$ be a fundamental domain for
${\G_j/\G_{U_{ij}}}$ in ${N/N_{U_{ij}}}$. By Lemma \ref{group
action} we have
\begin{align}
I &= \kappa \sum_{X \in \mathcal{C}_{\beta,i,j}}  \det g(X) \sum_{\g \in {\G_j/\G_{U_{ij}}}} \int_{\mathcal{D}''} \varphi_0(a(t)^{-1}n(b)^{-1}\g^{-1}X) db\\
  &= \kappa \sum_{X \in  \mathcal{C}_{\beta,i,j}}  \det g(X) \int_{   {N/N_{U_{ij}}}  } \varphi_0(a(t)^{-1}n(b)^{-1}X) db.
\end{align}

So we have proved
\begin{proposition}
\begin{multline*}
     A_{i,j}(t) =  2^{n/2} \det(v)^{\frac12} {(-1)^{(p-n)}}
     t^{p-n-1} \left( \int_{C_{U_{ij}}} \eta \right) \\
    \times \sum_{X \in \mathcal{C}_{\beta,i,j}}  \det g(X) \int_{N'}
    \varphi_0( iv, a(t)^{-1} n(b')^{-1}X) db'.
\end{multline*}
\end{proposition}

For the integral, we write
\begin{equation}
I(t,X) = \int_{N'}  \exp(-\pi tr(X,X)_{iv,Z(t,b')}) db',
\end{equation}
and it is not hard to see (\cite{KShin} Lemma 5.3)

\begin{lemma}\label{matrix}
\begin{equation*}
I(t,X) = 2^{-(n-1)/2} t^{n-1} e(\frac12 tr(i\beta v)) (\det
v)^{-\frac12} |\det g(X)|^{-1} \xi^{-\frac12}
\exp(-\frac{\pi}{2}t^{-2}\xi^{-1}),
\end{equation*}
where
\begin{equation*}
\xi= \xi(X) = \frac{\det v[^tg_1(X)]}{\det v \, \det g(X)^2}.
\end{equation*}
Here $ g(X) = \left(  \begin{smallmatrix} g_0(X) \\ g_1(X)
\end{smallmatrix} \right)$, where $g_0(X)$ is the first row and
$g_1(X)$ an $(n-1)$ by $n$ matrix.
\end{lemma}

\vspace{.5cm}

We are now ready to compute $\theta(\eta,\beta,U_{ij})$. We have
\begin{multline}\label{formula3}
    \theta(\eta,\beta,U_{ij})  = \frac{(-1)^{(n-1)}}{2^{1/2}}
     \left( \int_{C_{U_{ij}}}  \eta \right) e\left(\frac12 tr(i\beta v)\right) \\
        \times \sum_{X \in \mathcal{C}_{\beta,i,j}} \int_0^{\infty} sgn\det(g(X))
               \xi^{-\frac12}
               \exp\left(-\frac{\pi}{2} t^{-2} \xi^{-1} \right)  t^{-2} dt
\end{multline}
At this point interchanging of summation and integration in
(\ref{formula3}) is not allowed. Instead, we define for $ s \in
\C$,
\begin{equation}
I(s) = \int_0^{\infty}\sum_{X \in \mathcal{C}_{\beta,i,j}}
sgn\det(g(X)) \xi^{-\frac12}
          \exp\left(-\frac{\pi}{2} t^{-2} \xi^{-1} \right)  t^{-2-s} dt.
\end{equation}
Via a similar argument as in Prop. \ref{unfolding}, the sum is
rapidly decreasing as $t \to \infty$ so that $I(s)$ is entire and
for $Re(s) > 1$ we can interchange integration and summation by an
argument similar to Lemma \ref{estimate}(iv);  see also the proof
of Prop.\ref{Dirichlet} below. Noting $ sgn\det(g(X)) =
\epsilon(X)$ in the notation of Section 3, we obtain
\begin{align}
I(s) &= \sum_{X \in \mathcal{C}_{\beta,i,j}} \epsilon(X)
\xi^{-\frac12}
      \int_0^{\infty}
          \exp\left(-\frac{\pi}{2} t^{-2} \xi^{-1} \right)  t^{-2-s} dt \\
   &= 2^{(s-1)/2} \pi^{-(s+1)/2} \G\left(\frac{1+s}{2}\right)
     \sum_{X \in \mathcal{C}_{\beta,i,j}} \epsilon(X)\xi^{\frac{s}2}.
\end{align}

The above series is (up to the factor $\det(v)^{-s}$) the
Dirichlet series
\begin{equation}
\Omega(s,v,\beta) = \sum_{X \in \mathcal{C}_{\beta,i,j}}
f(s,v,g(X)),
\end{equation}
where $f: \C \times \mathbb{P}_n \times GL_n(\R) \longrightarrow
\C$ is given by \begin{equation} f(s,v,g) = \frac{sgn \det
g}{|\det g|^{s}} \det v[^tg_1]^{s}
\end{equation}
and $ g = \left( \begin{smallmatrix} g_0 \\ g_1 \end{smallmatrix}
\right)$ as in Lemma \ref{matrix}.

\begin{proposition}[\cite{KShin}]\label{Dirichlet}
 $\Omega(s,v,\beta)$ has an analytic continuation into the entire complex plane and
\begin{equation}
\Omega(0,v,\beta) =  - \sum_{X \in \mathcal{C}^{red}_{\beta,i,j}}
                             \mathbf{B}_1(\nu(X)) sgn\det(g(X)).
\end{equation}
\end{proposition}

\begin{proof}
Again we take $m \in SL_n(\Z)$ as in the proof of Lemma
\ref{cosetreps} and obtain
\begin{align}
\Omega(s,v,\beta) &= \sum_{X \in \mathcal{C}_{\beta,U_{ij},h}} f(s,v,g(X)) \\
                    &=  \sum_{Y \in \mathcal{C}_{\beta_0,U_{ij},k}} f(s,v',g(Y)).
\end{align}
in the notation of the proof of Lemma \ref{cosetreps} with $v' =
{^tm}^{-1}vm^{-1}$. Note $g(Y)= y_{01} \det(Y'_1)$.  Then
\begin{multline}
 \sum_{Y \in \mathcal{C}_{\beta_0,U_{ij},k}} f(s,v',g(Y)) =
 \sum_{y_{01} \equiv k_{01}} sgn(y_{01}) |y_{01}|^{-s}
   \sum_{Y'_1} sgn\det(Y'_1) \frac{\det v'[{^t Y'}_1]^s}{|\det Y_1'|^{s}}.
\end{multline}
The sum over $Y'_1$ is finite and can be evaluated for $s=0$,
while the first is equal to $H(k_{01},s) - H(1-k_{01},s)$, where
$H(x,s) = \sum_{n=0}^{\infty} (x+n)^{-s}$ is the Hurwitz
$\zeta$-function. The series converges for $Re(s) >1$ and has an
 analytic continuation to the whole complex plane. Observing
$H(x,0) = \tfrac12 - x= \mathbf{B}_1(x)$ for $x \in [0,1)$
finishes the proof of the proposition. Note here that the two
Hurwitz  $\zeta$-functions correspond to the two reduced elements
in one $\Z$-class in $\mathcal{C}_{\beta_0,U_{ij},k}$.
\end{proof}

Hence
\begin{equation}
I(0)=      - 2^{-1/2} \sum_{X \in \mathcal{C}^{red}_{\beta,i,j}}
                             \mathbf{B}_1(\nu(X)) sgn\det(g(X))
\end{equation}
and therefore
\begin{equation}
\theta(\eta,\beta,U_{ij})  = (-1)^{n}
                       \sum_{X \in \mathcal{C}^{red}_{\beta,i,j}}
                          \frac12 \mathbf{B}_1(\nu(X)) \epsilon(X)
                        e^{-2\pi  tr(\beta v)}.
\end{equation}
Considering (\ref{formula1}) in conjunction with the definition of
the cycle $C_{\beta}$ from Section \ref{cycles} concludes the
proof of Theorem \ref{semidef}!

\end{proof}

\section{The Positive Definite Fourier Coefficients}

\subsection{The defect $\delta_{\beta}(\eta)$}

For $\beta >0$, the main point of \cite{KMI,KMII,KMCan} (in much greater generality) is that the Fourier coefficient
\begin{align}\label{beta-coeff}
\theta_{\beta} &= 
\sum_{ X \in \Omega_{\beta} \cap (L^n+h)} \varphi(iv,Z,X) e^{-2\pi tr(\beta v)} \\
&= 
\sum_{ X \in \G \back \Omega_{\beta} \cap (L^n+h)} \sum_{\g \in \G_X \back \G} 
\g^{\ast} \varphi(iv,Z,X) e^{-2\pi tr(\beta v)}
\end{align}
is a Poincar\'e dual form for the composite cycle $C_{\beta}$, i.e.,
\begin{equation}\label{Dual}
a_{\beta}(\eta) = \int_{M} \eta \wedge \theta_{\beta} = \int_{C_{\beta}} \eta 
\end{equation}
for all $\eta \in {Z}_{rd}^{k}(M)$, the rapidly decreasing closed $k$-forms in $M$, and $k=p-n$.
(Actually, the case $n=p-1$ is not treated there, but for $p=2$ and $n=1$ we will show below that this is still true). 

Furthermore, by Stokes' Theorem, (\ref{Dual}) also holds on the space of relative coboundaries $B^{k}(\overline{M}, \partial \overline{M})$. Slightly more general we have
\begin{lemma}\label{PDLemma}
If $\eta$ is an exact $k$-form vanishing on $\partial \overline{M}$, then (\ref{Dual}) holds.
\end{lemma}

\begin{proof}
We consider the inclusion $i:M \longrightarrow \overline{M}$ and note that as a consequence of the relationship between duality on $H^{\ast}(\overline{M})$ and duality on  $H^{\ast} (\partial \overline{M})$, we obtain
\begin{equation}\label{PDboundary}
[i^{\ast} \theta_{\beta}] = PD[\partial_{\ast}[C_{\beta}]],
\end{equation}
see e.g. \cite{Bredon}, Th. 9.2, p. 357.
We write $\eta = d \omega$, whence the restriction of $\omega$ to $
 \partial \overline{M}$ is closed. Then
\begin{equation}\label{crucialpoint}
\int_M \eta \wedge \theta_{\beta} = \int_{\partial\overline{M}} \omega \wedge \theta_{\beta} = \int_{\partial C_{\beta}} \omega = \int_{C_{\beta}} \eta.
\end{equation}
Here the second equality follows from (\ref{PDboundary}) and that $\omega$ is closed on $\partial \overline{M}$. 
\end{proof}

However, (\ref{Dual}) will not hold for all $\eta \in {Z}^k(\overline{M})$ unless $C_{\beta}$ is compact (which can only happen for $k=p-n \leq 4$). 

In fact, we define the defect 
\begin{equation}
\delta_{\beta}(\eta) = a_{\beta}\left( \eta\right) - \int_{C_{\beta}} \eta
\end{equation}
for $\eta \in {Z}^{k}(\overline{M})$. By the above discussion, $\delta_{\beta}$ factors through 
\begin{equation}\label{spaces}
{Z}_{rd}^{k}(M) +
{B}^{k}(\overline{M}, \partial \overline{M})= 
{Z}^{k}(\overline{M}, \partial \overline{M}), 
\end{equation}
the closed forms vanishing at the boundary. The equality in (\ref{spaces}) follows from the fact that 
$H^{k}(\overline{M}, \partial \overline{M},\C) \simeq 
H_c^{k}({M},\C)$ has a system of representatives consisting of rapidly decreasing forms. So we proved
\begin{lemma}
For $n<p-1$ or $p=2$ and $n=1$, $\delta_{\beta}$ descends to a map
\begin{equation}\label{map}
\delta_{\beta}: \frac{Z^k(\overline{M})}
{Z^{k}(\overline{M}, \partial \overline{M})} \longrightarrow \C.
\end{equation}
\end{lemma}
We take a neighborhood $U$ of the boundary of $\overline{M}$ such that $\partial \overline{M}$ is a deformation retract of $U$ and obtain a projection map $\pi: U \longrightarrow  \partial \overline{M}$. We pick a smooth 'bump' function $\rho$ on $M$ supported in $U$ 
with $\rho \vert_V = 1$ for another neighborhood $V \subset U$ and define a map
\begin{equation}
\iota: A^{k-1}(\partial \overline{M}) \longrightarrow Z^k(\overline{M})
\end{equation}
by $\iota (\omega) = d \left(\rho \pi^{\ast}(\omega)\right)$ on $U$ and $\iota(\omega) = 0$ elsewhere. Note that 
 $\iota (\omega) \vert_{ \partial \overline{M}} = d \omega$. 
\begin{lemma}
We have the following exact sequence
\begin{equation} \label{exactseq}
0 \longrightarrow \frac{ A^{k-1}(\partial \overline{M})}{ Z^{k-1}(\partial \overline{M})} \stackrel{\bar{\iota}}{\longrightarrow} 
\frac{Z^k(\overline{M})}
{Z^{k}(\overline{M}, \partial \overline{M})}  \stackrel{\bar{r}}{\longrightarrow} H^k(\partial \overline{M},\C).
\end{equation}
Here $\bar{r}$ is the quotient of the restriction map $r: 
Z^k(\overline{M}) \longrightarrow   Z^k(\partial \overline{M})$ to the boundary; in general this is not surjective.
Also note that $\bar{\iota}$ is independent of the choices involved, so that (\ref{exactseq}) is intrinsic to the situation. 
\end{lemma}

We can therefore study the map (\ref{map}) via the exact sequence (\ref{exactseq}). 
\begin{proposition}
 $\delta_{\beta}$ is not identically zero on (the image of) $ \frac{ A^{k-1}(\partial \overline{M})}{ Z^{k-1}(\partial \overline{M})}$.
\end{proposition}

\begin{proof}
Let  $\omega \in A^{k-1}(\partial \overline{M})$. Then the calculation (\ref{crucialpoint}) for $\iota(\omega)$ is no longer valid (unless $\omega$ is closed) - and it is clear that there are examples so that (\ref{crucialpoint}) does not hold, i.e., $\delta_{\beta}(\iota(\omega)) \ne 0$.
For Riemann surfaces, we make this more explicit in the next section. 
\end{proof}

It is very tempting to investigate the other piece coming from $
 H^k(\partial \overline{M},\C)$ using Eisenstein cohomology.
We carry this out for the Riemann surface case.

\subsection{The defect for Riemann surfaces}

For the remainder of this section we consider the special case of
$SO_0(2,1)$. In particular, we prove the Theorems \ref{THE} and \ref{THF}.

There is a double covering $SL_2(\R) \longrightarrow SO_0(2,1)$
and the symmetric space $D$ is just the upper half plane $\h$. We
therefore work with $SL_2$ in this section. Accordingly, we change
notation and write $z=x+iy \in D \simeq \h$ for the orthogonal
variable. We write $dx_i$ for the basic differential form of the
boundary component $\mathbb{T}_i$ of the Borel-Serre
compactification corresponding to the cusp $\ell_i$. 
Hence $dx_i = (g_i^{-1})^{\ast}dx$.
Finally, for convenience we assume that $\G = \G(N) $, the principal congruence
subgroup. Hence the width of all cusps is equal to $N$.

\vspace{.5cm}

We first illustrate that the defect is not identically zero on
 $\frac{ A^{0}(\partial \overline{M})}{ Z^{0}(\partial \overline{M})}$.

By Theorem \ref{growth}, the restriction of $\theta(\tau,z)$ to a
boundary component $\mathbb{T}_i$ is given by
\begin{equation}
\theta(\tau,z)|_{\mathbb{T}_i} = 
\left(  \sum_{X \in W_i \cap L+h} P_2(X) e^{\pi i (X,X) \tau} \right) dx_i =: 
\theta_i(\tau) dx_i
\end{equation}
(For the isotropic line $\ell_i$ defining a cusp, $W_i
= \ell_i^{\perp}/\ell_i $ is one-dimensional, and
identifying $W_i(\R)$ with $\R$ we have $P_2(X) e^{\pi i (X,X)} = Xe^{\pi i X^2}$.)

Pick a function $f \in A^0(\partial{M}) = \oplus_i A^0(\mathbb{T}_i)$ only supported at the cusp $\ell_0$. Then
\begin{equation}
\int_M \iota(f) \wedge \theta(\tau,z) = \left( \int_{\mathbb{T}_0} f(x)dx \right) \theta_0(\tau);
\end{equation}
i.e., $A_{\beta}(\iota(f))$ is up to a factor the integral of $f$ over the whole boundary circle. On the other hand,
\begin{equation}
\int_{C_{\beta}} \iota(f)  = 
\int_{ \partial C_{\beta} \cap \mathbb{T}_0} f,
\end{equation}
which is the oriented sum of the evaluations of $f$ at the boundary points of $C_{\beta}$. It is certainly easy to find $f$, where these two terms are not the same; i.e., $\delta_{\beta}(\iota(f)) \neq 0$.

\vspace{.5cm}

We briefly review the relevant facts for Eisenstein  series and
Eisenstein cohomology needed.

We introduce the tangential Eisenstein series for the cusp
$\ell_i$,
\begin{equation}
E_i^T(s,z) = \sum_{\g \in \G_i \back \G}
Im(g_i^{-1}\g)^{\ast}(y^sdx)
\end{equation}
with $s \in \C$. We easily see
\begin{equation}\label{EisFormula}
E_i^T(s,z) = \frac{1}{2y} \left( E_i(s+1,z)_{-2}dz +
E_i(s+1,z)_2d\bar{z} \right)
\end{equation}
with
\begin{equation}
E_i(s,z)_n = \sum_{\g \in \G_i \back \G} Im(g_i^{-1}\g z)^s
\lambda(g_i^{-1}\g,z)^n,
\end{equation}
where $ \lambda(g,z)= \tfrac{cz+d}{|cz+d|}$ for $g = \left(
\begin{smallmatrix} a &b \\ c&d
\end{smallmatrix} \right)$.

The following theorem is well known, convenient references are
\cite{Kubota} and \cite{Freitag}.

\begin{theorem}
\begin{itemize}
\item[(i)]
The series $E_i(s,z)_{\pm 2}$ converge for $s >1$ and have a
meromorphic continuation to $\C$. At $s=1$, $E_i(s,z)_{\pm 2}$ are
holomorphic, and the Fourier expansions $E_{ij}(1,z)_{\pm 2}$ at a
cusp $\ell_j$ are given by
\begin{align}
\frac1{y} E_{ij}(1,z)_{- 2} &= \left(\delta_{ij} +
\frac1{y}a_{ij}(0)\right) + \sum_{m=1}^{\infty} a_{ij}(m) e^{2\pi
i mz/N }, \\
\frac1{y} E_{ij}(1,z)_{ 2} &= \left(\delta_{ij} +
\frac1{y}a_{ij}(0)\right) + \sum_{m=1}^{\infty}
\overline{a_{ij}(m)} e^{2\pi i m\bar{z}/N }.
\end{align}

\item[(ii)]
The tangential Eisenstein series $E_i^T(s,z)$ is holomorphic at
$s=0$ and defines a harmonic $1$-form on ${M}$, which extends to the boundary. For two different cusps $i$ and $j$, the difference
\begin{equation}
  E_i^T(0,z) -  E_j^T(0,z)
\end{equation}
is closed, and its restriction to the boundary is $dx_i - dx_j \in Z^1(\partial M) = \oplus_k Z^1(\mathbb{T}_k) $. We
call the space of all linear combination of tangential Eisenstein
series consisting of closed forms $\mathcal{E}_0$.

\item[(iii)]
The cohomology $H^1(M,\C)$ splits as
\[
H^1(M,\C) = H^1_{(2)}(M,\C) \oplus H^1_{Eis}(M,\C),
\]
where $H^1_{Eis}(M,\C)$ is the image of $\mathcal{E}_0$ in
$H^1(M,\C)$, while $ H^1_{(2)}(M,\C)$ is the $L_2$-cohomology. Its
classes can be represented by weight-$2$ cusp forms. Note $
H^1_{(2)}(M,\C) \simeq H^1_{!}(M,\C):= Im\left( H^1_c(M,\C)
\rightarrow H^1(M,\C) \right)$.

\end{itemize}
\end{theorem}

Theorem \ref{THE} now will follow from the vanishing of the defect $\delta_{\beta}$ for tangential Eisenstein  series and weight-$2$ cusp forms.
Via (\ref{beta-coeff}) we have to show
\begin{equation}\label{crucial}
\int_{\G \back B} \eta \wedge  \sum_{\g \in \G_X \back \G}  \g^{\ast}
\varphi(iv,z,X) = e^{-\pi (X,X)} \int_{C_X} \eta.
\end{equation}
for $(X,X) >0$. $X^{\perp}$ has signature $(1,1)$ and therefore 
the stabilizer $\G_X$ is either infinitely cyclic or trivial. In the first case, the cycle $C_X$ is a closed geodesic and (\ref{crucial}) holds for any $1$-form $\eta$. When the stabilizer is trivial, the cycle $C_X$ is an infinite geodesic joining two cusps. 
\begin{theorem}\label{ThomEisenstein}
Assume $C_X$ is an infinite geodesic. Then
\begin{equation}\label{Anna}
\int_{\G \back B} E_i^{T}(0,z) \wedge   \sum_{\g \in \G} \g^{\ast}
\varphi(z,X)= \ e^{-\pi (X,X)} \int_{C_X} E_i^T(0,z).
\end{equation}
\end{theorem}

\begin{proof}
First note that unfolding in (\ref{Anna}) is not allowed. Recall
we have a Witt basis $u_0,w,u_0'$ for $V$, and we can assume that
$X=2au_0 + bw$ with $a \in \Q$ and $b \in \Q_+$ so that $C_X$ is
the geodesic joining the cusps $\infty$ and $\tfrac{a}{b} \in \Q$.
The stabilizer of the cusp $\infty$ is $\G_{\infty} =
\G_{\infty}(N) = \left\{ \left(
\begin{smallmatrix}
1&k \\ 0&1
\end{smallmatrix}
\right)
: k \in N\Z \right\}$.

We have
\begin{equation}
\int_{\G \back B} E_i^{T}(0,z) \wedge   \sum_{\g \in \G} \g^{\ast}
\varphi(z,X)= \int_{\G \back B} E_i^{T}(0,z) \wedge \sum_{\g \in
\G_{\infty} \back \G} \sum_{k \in N\Z } \g^{\ast} \varphi(z,X+
2kbu_0).
\end{equation}
We introduce a holomorphic function $I(s)$ for $s \in \C$ by
\begin{equation}
I(s)= \int_{\G \back B} E_i^{T}(0,z) \wedge   \sum_{\g \in
\G_{\infty} \back \G} \sum_{k \in N\Z } \g^{\ast} \biggl( y^s
\varphi(z,X+2kbu_0) \biggr)
\end{equation}
and unfold
\begin{equation}
I(s)  = \int_{\G_{\infty} \back B} E_i^{T}(0,z) \wedge \sum_{k \in
N \Z }   y^s \varphi(z,X+2kbu_0).
\end{equation}
To justify this we first need some explicit formulae for $\varphi$. We have
\begin{align}
\varphi(z,X+2kbu_0)
&= \sqrt{2} b e^{-\pi (2\frac{(a-xb+kb)^2}{y^2} + b^2)} \frac{dx}y
+ \sqrt{2}\left(\frac{a-xb+kb}{y}\right)
e^{-\pi (2\frac{(a-xb+kb)^2}{y^2} + b^2)} \frac{dy}y \\
&=: \varphi_1(k,X)dx + \varphi_2(k,X)dy,
\end{align}
so that the Fourier transform with respect to $k$ is given by
\begin{equation}
\widehat{\varphi}(z,X+2kbu_0) =
\widehat{\varphi_1}(k,X)dx + \widehat{\varphi}_2(k,X)dy
\end{equation}
with
\begin{align}
\widehat{\varphi_1}(k,X)  &=
 e^{-\pi b^2} e^{-\pi \frac{(yk)^2}{2b^2}} e^{-2\pi i kx} e^{2\pi i \frac{a}{b} k}, \\
\widehat{\varphi}_2(k,X)  &=
 -i \frac{ky}{2b^2} e^{-\pi b^2}  e^{-\pi \frac{(yk)^2}{2b^2}} e^{-2\pi i kx} e^{2\pi i \frac{a}{b} k}.
\end{align}
By Poisson summation and (\ref{EisFormula}) we obtain
\begin{align}
I(s) &=  \int_{\G_{\infty} \back B} \frac{1}{2y} (E_i(1,z)_{-2} +
E_i(1,z)_2) \left( \frac{1}{N}  \sum_{k \in \frac1N \Z } \widehat{\varphi_2}(k,X) \right) y^s
 \\ &\qquad \qquad +
 \frac{-i}{2y} (E_i(1,z)_{-2} -
E_i(1,z)_2) \left( \frac1N \sum_{k \in \frac1N \Z } \widehat{\varphi_1}(k,X) \right) y^s dxdy,
\end{align}
and this is rapidly decreasing as $y \to \infty$ and of moderate
growth as $ y \to 0$. Hence unfolding is valid for $Re(s)$
sufficiently large. We pick the standard fundamental domain for
$\G_{\infty} \back B$ and integrate w.r.t. $x$. This picks out the
$0$-th Fourier coefficient:
\begin{align}
I(s) &=  \frac{-iy}{4b^2N} e^{-\pi b^2} \int_0^{\infty}
\sum_{m=1}^{\infty}  (\chi(m) a_i(m) -\overline{\chi(m)a_i(m)}) m
e^{-\pi \frac{(my)^2}{2(bN)^2} -2\pi
\frac{my}{N}} y^s dy \\
&\qquad \qquad - \frac{i}{2} e^{-\pi b^2} \int_0^{\infty}
\sum_{m=1}^{\infty} (\chi(m)a_i(m) -\overline{\chi(m)a_i(m)})
 e^{-\pi \frac{(my)^2}{2(bN)^2} -2\pi
\frac{my}{N}} y^s dy,
\end{align}
where $\chi(m) =  e^{2\pi i\frac{a}{b} m/N}$. Hence
\begin{align}
I(s) &=  \frac{-i e^{-\pi b^2/2}}{2 \sqrt{2}b} \sum_{m=1}^{\infty}
(\chi(m)a_i(m) -\overline{\chi(m)a_i(m)}) \\
& \hspace{2cm} \times\int_0^{\infty} \left( \frac{my}{\sqrt{2}bN}
+ \sqrt{2b} \right) e^{-\pi \left( \frac{my}{\sqrt{2}bN} +
\sqrt{2}b \right)^2} y^s dy \notag
\\
&= \frac{-i Ne^{-\pi b^2/2}}{4} (\sqrt{2}bN)^s
\left(L(E_i(1,z)_{-2},\chi,s+1) -L(E_i(1,z)_{2}, \overline{\chi},
s+1 ) \right) \\
& \hspace{2cm}  \times \int_{b^2/2}^{\infty} \left(  \sqrt{t} -
\sqrt{2}b \right)^s e^{-\pi t} dt, \notag
\end{align}
where $L(E_i(1,z)_{\pm 2},..., s )$
 are the (twisted) $L$-functions attached
to  $E_i(1,z)_{-2}$ and $E_i(1,z)_{2}$. Specializing to $s=0$ finally
gives
\begin{equation}
\int_{M} E_i^{T}(0,z) \wedge   \sum_{\g \in \G} \g^{\ast}
\varphi(Z,X) = \frac{-iN}{4\pi} (L(E_i(1,z)_{-2},\chi,1)
-L(E_i(1,z)_{2}, \overline{\chi}, 1 )) e^{-\pi b^2}.
\end{equation}
But now one easily checks that
\begin{equation}
\int_{C_X} E_i^{T}(0,z) = \frac{-iN}{4\pi}
(L(E_i(1,z)_{-2},\chi,1) -L(E_i(1,z)_{2}, \overline{\chi}, 1 )).
\end{equation}
This proves the theorem.
\end{proof}

\begin{remark}
The given proof (or a slightly simpler version of it) also works
for $\eta= f(z) dz$ with $f(z)$ a weight-$2$ cusp form. This is important,
 since for $C_X$ an infinite geodesic, the proof of
the basic identities (\ref{crucial}) and (\ref{Dual}) for $\eta$ rapidly decreasing
is actually not included in \cite{KMCan},\cite{KM90}.
\end{remark}

\vspace{.5cm}

Because of Stokes' theorem we have $ \int_{\partial{\overline{M}}}
\theta(\tau,z) =0$ and therefore $\sum_i \theta_i(\tau)=0$. Thus
\begin{equation}
Eis(\theta)(\tau,z) := \sum_i \theta_i(\tau) E^T_i(0,z)
\end{equation}
defines a closed differential form in $\overline{M}$ with values
in the holomorphic cusp forms of weight $3/2$, and we define the
truncated theta function
\begin{equation}
\theta^c(\tau,z) = \theta(\tau,z) - Eis(\theta,z).
\end{equation}
So $\theta^c(\tau)$ is per construction a rapidly decreasing
closed differential $1$-form in $M$ with values in the
non-holomorphic modular forms of weight $3/2$. We write
\begin{equation}
\theta^c(\tau,z) = \sum_{\beta} \theta^c_{\beta}(v,z) e^{2\pi i
\tau}
\end{equation}
for the Fourier expansion. For $\eta = f(z)dz$ with $f(z)$ a
weight-$2$ cusp form we still have
\begin{equation}\label{cuspforms}
\int_M \eta \wedge \theta^c_{\beta} =  \int_M \eta \wedge
\theta_{\beta} = \int_{C_{\beta}} \eta,
\end{equation}
as cusp forms are orthogonal to Eisenstein series. (By Theorem
\ref{ThomEisenstein}, (\ref{cuspforms}) actually also holds for
tangential Eisenstein since it is not too had to show that the
integral of the wedge of two tangential Eisenstein series
vanishes.) This justifies the
\begin{definition}
We define $C^c_{\beta}$ to be the homology class dual to the
Fourier coefficient $\theta^c_{\beta}$.
\end{definition}
$C_{\beta}$ does not depend on $v$ since (\ref{cuspforms}) and Th.
\ref{ThomEisenstein} show that $ \int_M \eta \wedge
\theta^c_{\beta}$ indeed does not depend on $v$.

This discussion proves Theorem \ref{THF}.

\section{The Theta Integral over Special Cycles}

We can also define a lift
\begin{equation}
\Lambda(\tau,C_U) = \int_{C_U} \theta_{\varphi_n}(\tau,Z),
\end{equation}
where $C_U$ is the special cycle coming from a positive definite
subspace $U$ of dimension $p-n$ in $V$.  Note that $C_U$ has
dimension $n$.

\vspace{.5cm}

We write $L_U = L \cap U$ and $L_{U^{\perp}} = L \cap U^{\perp}$
and obtain a decomposition
\begin{equation}\label{dec99}
L^n + h = \sum_{i=1}^s \left(L_U^n + h'_i\right) \; + \;
\left(L_{U^{\perp}}^n + h_i''\right)
\end{equation}
with $h_i' \in \left(L_u^{\#}\right)^n$ and $h_i'' \in
\left(L_{U^{\perp}}^{\#}\right)^n$. By $\theta_{C_U}(\tau,
L_{U^{\perp}} + h_i'')$ we denote the top degree theta integral
$\Lambda(\tau,1) =  \int_{C_U} \sum_{X \in  L_{U^{\perp}}^n +
h_i''} \varphi_n(\tau,X)$ for the hyperbolic space $C_U$. Note
that the top degree lift $\Lambda(\tau,1)$ was computed in
\cite{KShin}.

\begin{proposition}\label{intersection}
With the above notation, we have
\[
\Lambda(\tau,C_U) = \sum_{i=1}^s \vartheta(\tau,  L_U + h'_i)
\theta_{C_U}(\tau, L_{U^{\perp}} + h_i'')
\]
where  $\vartheta(\tau,  L_U + h'_i) = \sum_{X \in  L_U^n + h'_i}
e^{ \pi i tr((X,X)\tau)}$ is the standard theta series of degree
$n$ for the positive definite space $U$.
\end{proposition}

\begin{proof}
Using the explicit formula for $\varphi_n = \varphi_{n,V}$ from
Section 4 one easily checks that under the pullback $i^{\ast}_U :
\mathcal{A}^n(B) \longrightarrow  \mathcal{A}^n(B_U)$
\begin{equation}
i^{\ast}_U  \varphi_{n,V} =  \varphi_{0,U} \otimes
\varphi_{n,U^{\perp}},
\end{equation}
where $ \varphi_{0,U}$ is just the standard Gaussian for the
positive definite space $U$. From this the proposition easily
follows.
\end{proof}

\begin{theorem}(Theorem \ref{THF})
\[
 \Lambda(\tau,C_U) = \sum_{\beta>0} [C_U.C_{\beta}]_{tr} e^{2\pi i
 tr(\beta\tau)}\; + \; (-1)^n \sum_{\substack{\beta \geq 0 \\ rk(\beta) =
 n-1}}  [C_U.C^s_{\beta}] e^{2\pi i tr(\beta\tau)},
\]
where $[C_U.C_{\beta}]_{tr}$ is the transversal intersection
number of $C_U$ and $C_{\beta}$ (i.e., the sum of the transverse
intersections counted with multiplicities $+1$ and $-1$) and $
[C_U.C^s_{\beta}]$ is the evaluation of the cohomological
intersection product.
\end{theorem}

\begin{proof}

First assume for simplicity that in (\ref{dec99}) we have $s=1$
and write $h'=h'_1$ and $h'' = h''_1$. Let $\beta \in Sym_n(\Q)$
be positive definite. It is easy to see that  a $(p-n)$-cycle
$D_Y$ with $\tfrac12(Y,Y) =\beta$ intersects $D_U$ transversely if
and only if the orthogonal projection of $Y$ onto $U^{\perp}$
spans has rank $n$. From that we conclude that the transversal
intersection number $[C_U.C_{\beta}]_{tr}$ is given by
\begin{equation}\label{transversal}
[C_U.C_{\beta}]_{tr} = \sum_{\substack{ \alpha_1 \geq 0 \\
\alpha_2 >0 \\ \alpha_1 + \alpha_2 = \beta}} r(\alpha_1,U) \,
\text{deg}(C_{\alpha_2,C_U}),
\end{equation}
where $r(\alpha_1,U)$ is the representation number of $\alpha_1$
in $L^n_U + h'$ and
\begin{equation}{\label{degree}}
\text{deg}(C_{\alpha_2,C_U}) = \sum_{X \in \G_U \back
\Omega_{\alpha_2} \cap L^n_{U^{\perp}} +h''} \epsilon(X)
\end{equation}
is the (weighted) degree of the $0$-cycle $C_{\alpha_2}$ in the
space $C_U$ defined by $\alpha_2$. But the right hand side of
(\ref{transversal}) is exactly the $\beta$-th Fourier coefficient
of $\vartheta(\tau, L_U + h')$ times the positive definite part of
$\theta_{C_U}(\tau, L_{U^{\perp}} + h'')$, which is given by
(\ref{degree}), see \cite{KShin}.

The statement for the singular coefficients is clear as
Th. \ref{THD} shows that the $\beta$-th coefficient of $\theta$ is
the Poincar\'e dual of the absolute cycle $C^s_{\beta}$. But one
can also see in the same combinatorial way as above that the
Fourier coefficient attached to a semidefinite $\beta$ represents
the intersection numbers at the Borel-Serre boundary of the
singular cycle $C^s_{\beta}$ and the boundary of $C_U$.
\end{proof}


\begin{thebibliography}{99}

\bibitem{Borel}
A. Borel, \emph{Introduction aux groupes arithm\'etiques}, Hermann
1969.

\bibitem{Bredon}
G. Bredon, \emph{Topology and Geometry}, GTM 139, Springer, 1993


\bibitem{BS}
A. Borel and J.-P. Serre, \emph{Corners and arithmetic groups};
Comment. Math. Helv. \textbf{48} (1973), 436-491.

\bibitem{Cogdell}
J. Cogdell, \emph{The Weil representation and cycles on Picard
modular surfaces}, J. reine u. angew. Math. \textbf{357} (1985),
115-137.

\bibitem{Freitag}
E. Freitag, \emph{Hilbert Modular Forms}, Springer, 1990


\bibitem{Funke}
J. Funke, \emph{Heegner divisors and non-holomorphic modular
forms}, accepted by Compositio Math.

\bibitem{HZ}
F. Hirzebruch and D. Zagier, \emph{Intersection numbers of curves
on Hilbert modular surfaces and modular forms of Nebentypus}, Inv.
Math. \textbf{36} (1976), 57-113.

\bibitem{Kubota}
T. Kubota, \emph{Elementary Theory of Eisenstein Series}, Halsted
Press, 1973.


\bibitem{KudlaHS}
S. Kudla, \emph{Holomorphic Siegel Modular Forms Associated to
$SO(n,1)$}, Math. Ann. \textbf{256} (1981), 517-534.



\bibitem{KShin}
S. Kudla, \emph{ On the integrals of certain singular
theta-functions}, J. Fac. Sci. Univ. Tokyo \textbf{28} (1982), no.
3, 439-465.


\bibitem{KMI}
S. Kudla and J. Millson, \emph{The Theta Correspondence and
Harmonic Forms I}, Math. Ann. \textbf{274} (1986), 353-378.


\bibitem{KMII}
S. Kudla and J. Millson, \emph{The Theta Correspondence and
Harmonic Forms II}, Math. Ann. \textbf{277} (1987), 267-314.

\bibitem{KMCan}
S. Kudla and J. Millson, \emph{Tubes, cohomology with growth
conditions and application to the theta correspondence}, Canad. J.
Math. \textbf{40} (1988), 1-37.



\bibitem{KM90}
S. Kudla and J. Millson, \emph{Intersection numbers of cycles on
locally symmetric spaces and Fourier coefficients of holomorphic
modular forms in several complex variables}, IHES Pub. \textbf{71}
(1990), 121-172.

\bibitem {Shintani}
T. Shintani, \emph{On construction of holomorphic cusp forms of
half integral weight}, Nagoya Math. J., \textbf{58} (1975),
83-126.

\bibitem{Weil}
A. Weil, \emph{Sur certains groupes d'op\'erateurs unitaires},
Acta. Math. \textbf{111} (1964), 143-211.










\end{thebibliography}
\end{document}